\documentclass{amsart}
\usepackage{amsmath}
\usepackage{amssymb}
\usepackage{amsfonts}

\newtheorem{theorem}{Theorem}[section]
\newtheorem{corollary}[theorem]{Corollary}

\newtheorem{definition}[theorem]{Definition}
\newtheorem{proposition}[theorem]{Proposition}
\newtheorem{remark}[theorem]{Remark}
\numberwithin{theorem}{section}
\newtheorem{acknowledgement}{Acknowledgement}
\input{tcilatex}

\begin{document}
\title[Unbounded operator valued local positive maps ]{Factorization
properties for unbounded local positive maps }
\author{Maria Joi\c{t}a}
\address{Department of Mathematics, Faculty of Applied Sciences, University
Politehnica of Bucharest, 313 Spl. Independentei, 060042, Bucharest, Romania
and Simion Stoilow Institute of Mathematics of the Romanian Academy, P.O.
Box 1-764, 014700, Bucharest, Romania}
\email{mjoita@fmi.unibuc.ro and maria.joita@pub.ro}
\urladdr{http://sites.google.com/a/g.unibuc.ro/maria-joita/}
\thanks{This work was partially supported by a grant of the Ministry of
Research, Innovation and Digitization, CNCS/CCCDI--UEFISCDI, project number
PN-III-P4-ID-PCE-2020-0458, within PNCDI III}

\begin{abstract}
In this paper we present some factorization properties for unbounded local
positive maps. We show that an unbounded local positive map $\phi $ on the
minimal tensor product of the locally $C^{\ast }$-algebras $\mathcal{A}$ and 
$C^{\ast }(\mathcal{D}_{\mathcal{E}}),$ where $\mathcal{D}_{\mathcal{E}}$ is
a Fr\'{e}chet quantized domain, that is dominated by $\varphi \otimes $id is
of the forma $\psi \otimes $id, where $\psi $ is an unbounded local positive
map dominated by $\varphi $. As an application of this result, we show that
given a local positive map $\varphi :$ $\mathcal{A}\rightarrow $ $\mathcal{B}%
,$ the local positive map $\varphi \otimes $id$_{M_{n}\left( \mathbb{C}%
\right) }$ is local decomposable for some $n\geq 2$ if and only if $\varphi $
is a local\ $CP$-map. Also, we show that an unbounded local $CCP$-map\textit{%
\ }$\phi $ on the minimal tensor product of the unital locally $C^{\ast }$%
-algebras $\mathcal{A}$ and $\mathcal{B},$ that is dominated by $\varphi
\otimes \psi $ is of the forma $\varphi \otimes \widetilde{\psi }$, where $%
\widetilde{\psi }$ is an unbounded local $CCP$- map dominated by $\psi $,
whenever $\varphi $ is pure.
\end{abstract}

\maketitle

\section{$\protect\bigskip $Introduction}

Locally \ $C^{\ast }$-algebras are generalizations of $C^{\ast }$-algebras,
on which the topology instead of being given by a single $C^{\ast }$-norm is
defined by an upward directed family of $C^{\ast }$-seminorms. The concrete
models for locally $C^{\ast }$-algebras are $\ast $-algebras of unbounded
linear operators on a Hilbert space. In the literature, the locally $C^{\ast
}$-algebras are studied under different names like pro-$C^{\ast }$-algebras
(D. Voiculescu, N.C. Philips), $LMC^{\ast }$-algebras (G. Lassner, K. Schm%
\"{u}dgen), $b^{\ast }$-algebras (C. Apostol) and multinormed $C^{\ast }$%
-algebras (A. Dosiev). The term locally \ $C^{\ast }$-algebra is due to A.
Inoue \cite{I}.

A locally $C^{\ast }$-algebra is a complete Hausdorff complex topological $%
\ast $-algebra $\mathcal{A}$ whose topology is determined by a upward
filtered family $\{p_{\lambda }\}_{\lambda \in \Lambda }\ $of $C^{\ast }$%
-seminorms defined on $\mathcal{A}$. An element $a\in \mathcal{A}$ is
positive if there exists $b\in \mathcal{A}$ such that $a=b^{\ast }b$ and it
is local positive if there exist $c,d\in \mathcal{A}$ and $\lambda \in
\Lambda $ such that $a=c^{\ast }c+d$ with $p_{\lambda }\left( d\right) =0$.
Thus, the notion of (local) completely positive maps appeared naturally
while studying linear maps between locally $C^{\ast }$-algebras. The
structure of strictly continuous completely positive maps between locally $%
C^{\ast }$-algebras is described in \cite{J}. The same is done for strongly
bounded completely positive maps of order zero \cite{MJ2}. Dosiev \cite{D1}
proved a Stinespring type theorem for unbounded local completely positive
and local completely contractive maps on unital locally $C^{\ast }$%
-algebras. A Radon-Nikodym type theorem for such maps was proved by Bhat,
Ghatak and Kumar \cite{BGK}. In \cite{MJ1}, we obtained a structure theorem
for unbounded local completely positive maps of local order zero.

Bhat and Osaka \cite{BO} proved some factorization properties for bounded
positive maps on $C^{\ast }$-algebras. In this paper, we extend the results
of \cite{BO} to unbounded local positive maps on locally $C^{\ast }$%
-algebras.

The paper is organized as follows. In Section 2 we gather some basic facts
on locally $C^{\ast }$-algebras, concrete models for locally $C^{\ast }$%
-algebras and unbounded local completely positive and local completely
contractive maps needed for understanding the main results of this paper. In
Section 3, we show that a linear map between locally $C^{\ast }$-algebras is
local $CP$ (completely positive) if and only if it is continuous and
completely positive (Proposition \ref{5}). Therefore, the local $CP$-maps%
\textit{\ }on unital locally $C^{\ast }$-algebras are exactly strictly
continuous completely positive maps on unital locally $C^{\ast }$-algebras 
\cite[Remark 4.4]{J} and the structure theorem \cite[Theorem 4.6]{J} is
valid for local\textit{\ }$CP$-maps\textit{\ }on unital locally $C^{\ast }$%
-algebras. As in the case of bounded completely positive maps on $C^{\ast }$%
-algebras, we show that two unbounded local $CCP$ (local completely
contractive and local completely positive)-maps on unital locally $C^{\ast }$%
-algebras determine an unbounded local $CCP$-map on the minimal tensor
product and a minimal Stinespring dilation for the tensor product map can be
obtained in terms of the minimal Stinespring dilations for each map\
(Proposition \ref{3}). In section 4, we show that if $\mathcal{D}_{\mathcal{E%
}}$ is a Fr\'{e}chet quantized domain, then an unbounded local positive map $%
\phi $ on the minimal tensor product of the locally $C^{\ast }$-algebras $%
\mathcal{A}$ and $C^{\ast }(\mathcal{D}_{\mathcal{E}}),$ that is dominated
by $\varphi \otimes $id$_{C^{\ast }(\mathcal{D}_{\mathcal{E}})}$ factorizes
as $\psi \otimes $id$_{C^{\ast }(\mathcal{D}_{\mathcal{E}})}$, where $\psi $
is an unbounded local positive map dominated by $\varphi \ $(Theorem \ref{6}%
). As an application of this result, we show that given a local positive map 
$\varphi :$ $\mathcal{A}\rightarrow $ $\mathcal{B}$, the local positive map $%
\varphi \otimes $id$_{M_{n}\left( \mathbb{C}\right) }$ is local decomposable
for $n\geq 2$ if and only if $\varphi $ is a local\ $CP$-map\ (Theorem \ref%
{7}). Also, we show that given unbounded\textbf{\ }local\textbf{\ }$CCP$%
-maps $\phi :\mathcal{A\otimes B}\rightarrow C^{\ast }(\mathcal{D}_{\mathcal{%
E\otimes F}}),$ $\varphi :\mathcal{A}\rightarrow C^{\ast }(\mathcal{D}_{%
\mathcal{E}})$ and $\psi :\mathcal{B}\rightarrow C^{\ast }(\mathcal{D}_{%
\mathcal{F}})$, if $\phi $ is dominated by $\varphi \otimes \psi $ and $%
\varphi $ is pure, then $\phi $ factorizes as $\varphi \otimes \widetilde{%
\psi }$, where $\widetilde{\psi }$ is an unbounded \ local $CCP$-map
dominated by $\psi $ (Theorem \ref{8}).

\section{Preliminaries}

Let $\mathcal{A}$ be a locally $C^{\ast }$-algebra with the topology defined
by the family of $C^{\ast }$-seminorms $\left\{ p_{\lambda }\right\}
_{\lambda \in \Lambda }.$

An element $a\in \mathcal{A}$ is \textit{bounded} if $\sup \{p_{\lambda
}\left( a\right) ;\lambda \in \Lambda \}<\infty $. The subset $b\left( 
\mathcal{A}\right) =\{a\in \mathcal{A};\left\Vert a\right\Vert _{\infty
}:=\sup \{p_{\lambda }\left( a\right) ;\lambda \in \Lambda \}<\infty \}\ $is
a $C^{\ast }$-algebra with respect to the $C^{\ast }$-norm $\left\Vert \cdot
\right\Vert _{\infty }.$ Moreover, $b\left( \mathcal{A}\right) $ is dense in 
$\mathcal{A}.$

Let us observe that $\mathcal{A}$ can be realized as a projective limit of
an inverse family of $C^{\ast }$-algebras as follows: For each $\lambda \in
\Lambda $, let $\mathcal{I}_{\lambda }=\{a\in \mathcal{A};p_{\lambda }\left(
a\right) =0\}$. Clearly, $\mathcal{I}_{\lambda }$ is a closed two sided $%
\ast $-ideal in $\mathcal{A}$ and $\mathcal{A}_{\lambda }=\mathcal{A}/%
\mathcal{I}_{\lambda }$ is a $C^{\ast }$-algebra with respect to the norm
induced by $p_{\lambda }$. The canonical quotient $\ast $-morphism from $%
\mathcal{A\ }$to $\mathcal{A}_{\lambda }$ is denoted by $\pi _{\lambda }^{%
\mathcal{A}}$. For each $\lambda _{1},\lambda _{2}\in \Lambda $ with $%
\lambda _{1}\leq \lambda _{2}$, there is a canonical surjective $\ast $%
-morphism $\pi _{\lambda _{2}\lambda _{1}}^{\mathcal{A}}:$ $\mathcal{A}%
_{\lambda _{2}}\rightarrow \mathcal{A}_{\lambda _{1}}$ defined by $\pi
_{\lambda _{2}\lambda _{1}}^{\mathcal{A}}\left( a+\mathcal{I}_{\lambda
_{2}}\right) =a+\mathcal{I}_{\lambda _{1}}$ for $a\in \mathcal{A}$. Then, $\{%
\mathcal{A}_{\lambda },\pi _{\lambda _{2}\lambda _{1}}^{\mathcal{A}}\}$\
forms an inverse system of $C^{\ast }$-algebras, since $\pi _{\lambda _{1}}^{%
\mathcal{A}}=$ $\pi _{\lambda _{2}\lambda _{1}}^{\mathcal{A}}\circ \pi
_{\lambda _{2}}^{\mathcal{A}}$ whenever $\lambda _{1}\leq \lambda _{2}$. The
projective limit%
\begin{equation*}
\lim\limits_{\underset{\lambda }{\leftarrow }}\mathcal{A}_{\lambda
}=\{\left( a_{\lambda }\right) _{\lambda \in \Lambda }\in
\tprod\limits_{\lambda \in \Lambda }\mathcal{A}_{\lambda };\pi _{\lambda
_{2}\lambda _{1}}^{\mathcal{A}}\left( a_{\lambda _{2}}\right) =a_{\lambda
_{1}}\text{ whenever }\lambda _{1}\leq \lambda _{2},\lambda _{1},\lambda
_{2}\in \Lambda \}
\end{equation*}%
of the inverse system of $C^{\ast }$-algebras $\{\mathcal{A}_{\lambda },\pi
_{\lambda _{2}\lambda _{1}}^{\mathcal{A}}\}$ is a locally $C^{\ast }$%
-algebra that can be identified with $\mathcal{A}$ via the map $a\mapsto
\left( \pi _{\lambda }^{\mathcal{A}}\left( a\right) \right) _{\lambda \in
\Lambda }.$

An element $a\in \mathcal{A}$ is \textit{self-adjoint} if $a^{\ast }=a$ and
it is\textit{\ positive} if $a=b^{\ast }b$ for some $b\in \mathcal{A}$. An
element $a\in \mathcal{A}$ is called \textit{local self-adjoint} if $%
a=a^{\ast }+c$, where $c\in \mathcal{A}$ such that $p_{\lambda }\left(
c\right) =0$ for some $\lambda \in \Lambda ,$ and \textit{local positive} if 
$a=b^{\ast }b+c$ where $b,c\in $ $\mathcal{A}$ such that $p_{\lambda }\left(
c\right) =0\ $ for some $\lambda \in \Lambda $. In the first case, we call $%
a $ as $\lambda $-self-adjoint, and in the second case, we call $a$ as $%
\lambda $-positive and write $a\geq _{\lambda }0$. We write, $a=_{\lambda }0$
whenever $p_{\lambda }\left( a\right) =0$. Note that $a\in \mathcal{A}$ is
local self-adjoint if and only if there is $\lambda \in \Lambda $ such that $%
\pi _{\lambda }^{\mathcal{A}}\left( a\right) $ is self adjoint in $\mathcal{A%
}_{\lambda }$ and $a\in \mathcal{A}$ is local positive if and only if there
is $\lambda \in \Lambda $ such that $\pi _{\lambda }^{\mathcal{A}}\left(
a\right) $ is positive in $\mathcal{A}_{\lambda }.$

Let $\mathcal{A}$ and $\mathcal{B}$ be two locally $C^{\ast }$-algebras with
the topology defined by the family of $C^{\ast }$-seminorms $\left\{
p_{\lambda }\right\} _{\lambda \in \Lambda }$ and $\left\{ q_{\delta
}\right\} _{\delta \in \Delta }$, respectively. For each $n\in \mathbb{N},$ $%
M_{n}(\mathcal{A})$ denotes the collection of all matrices of size $n$ with
elements in $\mathcal{A}$. Note that $M_{n}(\mathcal{A})$ is a locally $%
C^{\ast }$-algebra where the associated family of $C^{\ast }$-seminorms is
denoted by $\{p_{\lambda }^{n}\}_{\lambda \in \Lambda }.$

For each $n\in \mathbb{N}$, the $n$-amplification of a linear map $\varphi :%
\mathcal{A}\rightarrow \mathcal{B}$ is the map $\varphi ^{\left( n\right)
}:M_{n}(\mathcal{A})$ $\rightarrow $ $M_{n}(\mathcal{B})$ defined by 
\begin{equation*}
\varphi ^{\left( n\right) }\left( \left[ a_{ij}\right] _{i,j=1}^{n}\right) =%
\left[ \varphi \left( a_{ij}\right) \right] _{i,j=1}^{n}
\end{equation*}%
for $\left[ a_{ij}\right] _{i,j=1}^{n}\in M_{n}(\mathcal{A})$ .

A linear map $\varphi :\mathcal{A}\rightarrow \mathcal{B}$ is called :

\begin{enumerate}
\item \textit{local contractive} if for each $\delta \in \Delta $, there
exists $\lambda \in \Lambda $ such that%
\begin{equation*}
q_{\delta }\left( \varphi \left( a\right) \right) \leq p_{\lambda }\left(
a\right) \text{ for all }a\in \mathcal{A};
\end{equation*}

\item \textit{local positive} if for each $\delta \in \Delta ,$there exists $%
\lambda \in \Lambda $\ such that $\varphi \left( a\right) \geq _{\delta }0$
whenever $a\geq _{\lambda }0$ and $\varphi \left( a\right) =_{\delta }0$\ \
\ whenever $a=_{\lambda }0;$

\item \textit{local completely contractive }(\textit{local }$CC$\textit{-map}%
)\textit{\ }if for each $\delta \in \Delta $, there exists $\lambda \in
\Lambda $ such that 
\begin{equation*}
q_{\delta }^{n}\left( \varphi ^{\left( n\right) }\left( \left[ a_{ij}\right]
_{i,j=1}^{n}\right) \right) \leq p_{\lambda }^{n}\left( \left[ a_{ij}\right]
_{i,j=1}^{n}\right) \text{ }
\end{equation*}%
for all $\left[ a_{ij}\right] _{i,j=1}^{n}\in M_{n}(\mathcal{A})\ $and for
all $n\in \mathbb{N};$

\item \textit{local completely positive }(\textit{local }$CP$\textit{-map})%
\textit{\ }if for each $\delta \in \Delta $, there exists $\lambda \in
\Lambda $ such that $\varphi ^{\left( n\right) }\left( \left[ a_{ij}\right]
_{i,j=1}^{n}\right) \geq _{\delta }0\ $whenever $\left[ a_{ij}\right]
_{i,j=1}^{n}\geq _{\lambda }0$ and $\varphi ^{\left( n\right) }\left( \left[
a_{ij}\right] _{i,j=1}^{n}\right) =_{\delta }0\ \ $whenever $\left[ a_{ij}%
\right] _{i,j=1}^{n}=_{\lambda }0,\ $for all $n\in \mathbb{N}.$
\end{enumerate}

Throughout the paper, $\mathcal{H}$ is a complex Hilbert space and $B(%
\mathcal{H})$ is the algebra of all bounded linear operators on $\mathcal{H}$%
.

Let $(\Lambda ,\leq )$ be a directed poset. A \textit{quantized domain }in a
Hilbert space $\mathcal{H}$ is a triple $\{\mathcal{H};\mathcal{E};\mathcal{D%
}_{\mathcal{E}}\}$, where $\mathcal{E}=\{\mathcal{H}_{\lambda };\lambda \in
\Lambda \}$ is an upward filtered family of closed subspaces with dense
union $\mathcal{D}_{\mathcal{E}}=\tbigcup\limits_{\lambda \in \Lambda }%
\mathcal{H}_{\lambda }$ in $\mathcal{H\ }$\cite{D1}.

A quantized family $\mathcal{E}=\{\mathcal{H}_{\lambda };\lambda \in \Lambda
\}$ determines an upward filtered family $\{P_{\lambda };\lambda \in \Lambda
\}$ of projections in $B(\mathcal{H})$, where $P_{\lambda }$ is a projection
onto $\mathcal{H}_{\lambda }$.

We say that a quantized domain $\mathcal{F}=\{\mathcal{K}_{\lambda };\lambda
\in \Lambda \}$ \ of $\mathcal{H}$ with its union space $\mathcal{D}_{%
\mathcal{F}}$ and $\mathcal{K}$ $=\overline{\tbigcup\limits_{\lambda \in
\Lambda }\mathcal{K}_{\lambda }}$ is a \textit{quantized subdomian} of $%
\mathcal{E}$, if $\mathcal{K}_{\lambda }\subseteq \mathcal{H}_{\lambda }$
for all $\lambda \in \Lambda $. In this case, we write $\mathcal{F}\subseteq 
$ $\mathcal{E}$.

Let $\mathcal{E}^{i}=\{\mathcal{H}_{\lambda }^{i};\lambda \in \Lambda \}$ be
a quantized domain in a Hilbert space $\mathcal{H}^{i}$ for $i=1,2.$Given a
linear operator $V:\mathcal{D}_{\mathcal{E}^{1}}\rightarrow \mathcal{H}^{2}$%
, we write $V(\mathcal{E}^{1})\subseteq \mathcal{E}^{2}$ if $V(\mathcal{H}%
_{\lambda }^{1})\subseteq \mathcal{H}_{\lambda }^{2}$\ for all $\lambda \in
\Lambda $.

Let $\mathcal{H}$ and $\mathcal{K}$ be Hilbert spaces. A linear operator $T$ 
$:$ dom$(T)$ $\subseteq $ $\mathcal{H}$ $\rightarrow $ $\mathcal{K}$ is said
to be densely defined if dom$(T)$ is a dense subspace of $\mathcal{H}$. The
adjoint of $T$ is a linear map $T^{\bigstar }:$ dom$(T^{\bigstar })$ $%
\subseteq $ $\mathcal{K}\rightarrow \mathcal{H},$ where%
\begin{equation*}
\text{dom}(T^{\bigstar })=\{\xi \in \mathcal{K};\eta \rightarrow
\left\langle T\eta ,\xi \right\rangle _{\mathcal{K}}\ \text{is continuous
for every }\eta \in \text{dom}(T)\}
\end{equation*}
satisfying $\left\langle T\eta ,\xi \right\rangle _{\mathcal{K}%
}=\left\langle \eta ,T^{\bigstar }\xi \right\rangle _{\mathcal{H}}$ for all $%
\xi \in $dom$(T^{\bigstar })$ and $\eta \in $dom$(T).$

Let $\mathcal{E}=\{\mathcal{H}_{\lambda };\lambda \in \Lambda \}$ be a
quantized domain in a Hilbert space $\mathcal{H}$ and 
\begin{equation*}
C(\mathcal{D}_{\mathcal{E}})=\{T\in \mathcal{L}(\mathcal{D}_{\mathcal{E}%
});TP_{\lambda }=P_{\lambda }TP_{\lambda }\in B(\mathcal{H})\text{ for all }%
\lambda \in \Lambda \}
\end{equation*}%
where $\mathcal{L}(\mathcal{D}_{\mathcal{E}})$ is the collection of all
linear operators on $\mathcal{D}_{\mathcal{E}}$. If $T\in \mathcal{L}(%
\mathcal{D}_{\mathcal{E}})$, then $T\in C(\mathcal{D}_{\mathcal{E}})$ if and
only if $T(\mathcal{H}_{\lambda })\subseteq \mathcal{H}_{\mathcal{\lambda }}$
and $\left. T\right\vert _{\mathcal{H}_{\lambda }}\in B(\mathcal{H}_{\lambda
})$ for all $\lambda \in \Lambda $, and so $C(\mathcal{D}_{\mathcal{E}})$ is
an algebra. Let 
\begin{equation*}
C^{\ast }(\mathcal{D}_{\mathcal{E}})=\{T\in C(\mathcal{D}_{\mathcal{E}%
});P_{\lambda }T\subseteq TP_{\lambda }\text{ for all }\lambda \in \Lambda
\}.
\end{equation*}
If $T\in C(\mathcal{D}_{\mathcal{E}})$, then $T\in C^{\ast }(\mathcal{D}_{%
\mathcal{E}})$ if and only if $T(\mathcal{H}_{\lambda }^{\bot }\cap \mathcal{%
D}_{\mathcal{E}})\subseteq \mathcal{H}_{\lambda }^{\bot }\cap \mathcal{D}_{%
\mathcal{E}}$ for all $\lambda \in \Lambda .\ $

If $T\in C^{\ast }(\mathcal{D}_{\mathcal{E}})$, then $\mathcal{D}_{\mathcal{E%
}}$ $\subseteq $ dom$(T^{\bigstar })$. Moreover, $T^{\bigstar }(\mathcal{H}%
_{\lambda })\subseteq \mathcal{H}_{\lambda }$ for all $\lambda \in \Lambda $%
. Now, let $T^{\ast }=\left. T^{\bigstar }\right\vert _{\mathcal{D}_{%
\mathcal{E}}}$. It is easy to check that $T^{\ast }\in C^{\ast }(\mathcal{D}%
_{\mathcal{E}})$, and so $C^{\ast }(\mathcal{D}_{\mathcal{E}})$ is a unital $%
\ast $-algebra.

For each $\lambda \in \Lambda $, the map $\left\Vert \cdot \right\Vert
_{\lambda }:C^{\ast }(\mathcal{D}_{\mathcal{E}})\rightarrow \lbrack 0,\infty
)$, 
\begin{equation*}
\left\Vert T\right\Vert _{\lambda }=\left\Vert \left. T\right\vert _{%
\mathcal{H}_{\lambda }}\right\Vert =\sup \{\left\Vert T\left( \xi \right)
\right\Vert ;\xi \in \mathcal{H}_{\lambda },\left\Vert \xi \right\Vert \leq
1\}
\end{equation*}%
is a $C^{\ast }$-seminorm on $C^{\ast }(\mathcal{D}_{\mathcal{E}})$.
Moreover, $C^{\ast }(\mathcal{D}_{\mathcal{E}})$ is a locally $C^{\ast }$%
-algebra with respect to the topology determined by the family of $C^{\ast }$%
-seminorms $\{\left\Vert \cdot \right\Vert _{\lambda }\}_{\lambda \in
\Lambda }$ and $b(C^{\ast }(\mathcal{D}_{\mathcal{E}}))$ is identified with
the $C^{\ast }$-algebra $\{T\in B\left( \mathcal{H}\right) ;P_{\lambda
}T=TP_{\lambda }\ $for all $\lambda \in \Lambda \}$\ via\ the map $T\mapsto 
\widetilde{T}$, where $\widetilde{T}$ is the extension of $T$ to $\mathcal{H}
$ (see the proof of \cite[Lemma 3.1]{D})\textbf{.}

If $\mathcal{E}=\{\mathcal{H}_{\lambda };\lambda \in \Lambda \}$ is a
quantized domain in a Hilbert space $\mathcal{H}$, then $\mathcal{D}_{%
\mathcal{E}}$ can be regarded as a strict inductive limit of the direct
family of Hilbert spaces $\mathcal{E}=\{\mathcal{H}_{\lambda };\lambda \in
\Lambda \},$ $\mathcal{D}_{\mathcal{E}}=\lim\limits_{\rightarrow }\mathcal{H}%
_{\lambda }$, and it is called a locally Hilbert space (see \cite{I}). If $%
T\in C^{\ast }(\mathcal{D}_{\mathcal{E}})$, then $T$ is continuous \cite[%
Lemma 5.2 ]{I}.

For every locally $C^{\ast }$-algebra $\mathcal{A}$ there is a quantized
domain $\mathcal{E}$ in a Hilbert space $\mathcal{H}$ and a local isometric $%
\ast $-homomorphism $\pi :\mathcal{A\rightarrow }C^{\ast }(\mathcal{D}_{%
\mathcal{E}})$ \cite[Theorem 7.2]{D1}. This result can be regarded as an
unbounded analog of the Ghelfand-Naimark theorem.

\section{Local completely positive maps}

Let $\mathcal{A}$ and $\mathcal{B}$ be two locally $C^{\ast }$-algebras with
the topology defined by the family of $C^{\ast }$-seminorms $\left\{
p_{\lambda }\right\} _{\lambda \in \Lambda }$ and $\left\{ q_{\delta
}\right\} _{\delta \in \Delta }$, respectively.

\begin{proposition}
\label{4} Let $\varphi :\mathcal{A}\rightarrow $ $\mathcal{B}$ be a linear
map. If $\varphi $ is local positive, then $\varphi $ is continuous.
\end{proposition}

\begin{proof}
Since $\varphi $ is local positive, for each $\delta \in \Delta $, there
exists $\lambda \in \Lambda $ such that $\varphi \left( a\right) \geq
_{\delta }0\ $whenever $a\geq _{\lambda }0$ and $\varphi \left( a\right)
=_{\delta }0\ \ \ $whenever $a=_{\lambda }0$. Define the map $\varphi
_{\delta }^{+}:\left( \mathcal{A}_{\lambda }\right) _{+}\rightarrow \mathcal{%
B}_{\delta }$ by $\varphi _{\delta }^{+}\left( \pi _{\lambda }^{\mathcal{A}%
}\left( a\right) \right) =\pi _{\delta }^{\mathcal{B}}\left( \varphi \left(
a\right) \right) $, and extend it to a linear map $\varphi _{\delta }:%
\mathcal{A}_{\lambda }\rightarrow \mathcal{B}_{\delta }$. This map is
positive, and so continuous. Therefore, there is $C_{\delta }>0\ $such that%
\begin{equation*}
q_{\delta }\left( \varphi \left( a\right) \right) =\left\Vert \pi _{\delta
}^{\mathcal{B}}\left( \varphi \left( a\right) \right) \right\Vert _{\mathcal{%
B}_{\delta }}=\left\Vert \varphi _{\delta }\left( \pi _{\lambda }^{\mathcal{A%
}}\left( a\right) \right) \right\Vert \leq C_{\delta }\left\Vert \pi
_{\lambda }^{\mathcal{A}}\left( a\right) \right\Vert _{\mathcal{A}_{\lambda
}}=C_{\delta }p_{\lambda }\left( a\right)
\end{equation*}%
for all $a\in \mathcal{A}$.
\end{proof}

Proposition \ref{4} it is a particular case of \cite[Lemma 4.4]{D1}.

\begin{remark}
\label{1}Let $\varphi :\mathcal{A}\rightarrow $ $\mathcal{B}$ be a local
positive linear map. Then, for each $\delta \in \Delta $, there exist $%
\lambda \in \Lambda $ and a positive map $\varphi _{\delta }:\mathcal{A}%
_{\lambda }\rightarrow \mathcal{B}_{\delta }$ such that $\varphi _{\delta
}\left( \pi _{\lambda }^{\mathcal{A}}\left( a\right) \right) =\pi _{\delta
}^{\mathcal{B}}\left( \varphi \left( a\right) \right) $ for all $a\in 
\mathcal{A}$.
\end{remark}

\begin{proposition}
\label{5}Let $\mathcal{A}$ and $\mathcal{B}$\ be two locally $C^{\ast }$%
-algebras and $\varphi :\mathcal{A}\rightarrow $ $\mathcal{B}$ be a linear
map. Then $\varphi $ is local completely positive if and only if $\varphi $
is continuous and completely positive.
\end{proposition}

\begin{proof}
If $\varphi $ is local completely positive, then, by Proposition \ref{4},
for each $n\in \mathbb{N}$, the map $\varphi ^{\left( n\right) }$ is
continuous and by \cite[Proposition 2.1]{MJ1}, it is positive. Therefore, $%
\varphi $ is continuous and completely positive.

Conversely, suppose that $\varphi $ is continuous and completely positive.
Since $\varphi $ is continuous, for each $\delta \in \Delta $, there exist $%
\lambda \in \Lambda $ and a positive map $\varphi _{\delta }:\mathcal{A}%
_{\lambda }\rightarrow \mathcal{B}_{\delta }$ such that $\varphi _{\delta
}\left( \pi _{\lambda }^{\mathcal{A}}\left( a\right) \right) =\pi _{\delta
}^{\mathcal{B}}\left( \varphi \left( a\right) \right) $ for all $a\in 
\mathcal{A}$. Let $\left[ a_{ij}\right] _{i,j=1}^{n}\in M_{n}\left( \mathcal{%
A}\right) $ such that $\left[ a_{ij}\right] _{i,j=1}^{n}\geq _{\lambda }0$.
This means that there exist $\left[ b_{ij}\right] _{i,j=1}^{n},$ $\left[
c_{ij}\right] _{i,j=1}^{n}\in M_{n}\left( \mathcal{A}\right) $ such that $%
\left[ a_{ij}\right] _{i,j=1}^{n}=\left( \left[ b_{ij}\right]
_{i,j=1}^{n}\right) ^{\ast }\left[ b_{ij}\right] _{i,j=1}^{n}+$ $\left[
c_{ij}\right] _{i,j=1}^{n}\ $and $p_{\lambda }^{n}\left( \left[ c_{ij}\right]
_{i,j=1}^{n}\right) =0$. Then%
\begin{eqnarray*}
\pi _{\delta }^{M_{n}(\mathcal{B)}}\left( \varphi ^{\left( n\right) }\left( %
\left[ c_{ij}\right] _{i,j=1}^{n}\right) \right) &=&\left[ \pi _{\delta }^{%
\mathcal{B}}\left( \varphi \left( c_{ij}\right) \right) \right] _{i,j=1}^{n}=%
\left[ \varphi _{\delta }\left( \pi _{\lambda }^{\mathcal{A}}\left(
c_{ij}\right) \right) \right] _{i,j=1}^{n} \\
&=&\varphi _{\delta }^{\left( n\right) }\left( \pi _{\lambda }^{M_{n}\left( 
\mathcal{A}\right) }\left( \left[ c_{ij}\right] _{i,j=1}^{n}\right) \right) =%
\left[ 0\right] _{i,j=1}^{n}
\end{eqnarray*}%
and, since $\varphi \ $is completely positive,%
\begin{equation*}
\pi _{\delta }^{M_{n}(\mathcal{B)}}\left( \varphi ^{\left( n\right) }\left( %
\left[ a_{ij}\right] _{i,j}^{n}\right) \right) =\pi _{\delta }^{M_{n}(%
\mathcal{B)}}\left( \varphi ^{\left( n\right) }\left( \left( \left[ b_{ij}%
\right] _{i,j=1}^{n}\right) ^{\ast }\left[ b_{ij}\right] _{i,j=1}^{n}\right)
\right) \geq 0.
\end{equation*}%
If $\left[ a_{ij}\right] _{i,j}^{n}=_{\lambda }\left[ 0\right] _{i,j=1}^{n}$,%
$\ $then $p_{\lambda }^{n}\left( \left[ a_{ij}\right] _{i,j}^{n}\right) =0\ $%
and%
\begin{equation*}
\pi _{\delta }^{M_{n}(\mathcal{B)}}\left( \varphi ^{\left( n\right) }\left( %
\left[ a_{ij}\right] _{i,j=1}^{n}\right) \right) =\varphi _{\delta }^{\left(
n\right) }\left( \pi _{\lambda }^{M_{n}\left( \mathcal{A}\right) }\left( %
\left[ a_{ij}\right] _{i,j=1}^{n}\right) \right) =\left[ 0\right]
_{i,j=1}^{n}.
\end{equation*}%
Therefore, 
\begin{equation*}
\varphi ^{\left( n\right) }\left( \left[ a_{ij}\right] _{i,j}^{n}\right)
\geq _{\delta }0
\end{equation*}%
for all $\left[ a_{ij}\right] _{i,j=1}^{n}\in M_{n}\left( \mathcal{A}\right) 
$ such that $\left[ a_{ij}\right] _{i,j=1}^{n}\geq _{\lambda }0,\ $and $%
\varphi ^{\left( n\right) }\left( \left[ a_{ij}\right] _{i,j}^{n}\right)
=_{\delta }\left[ 0\right] _{i,j=1}^{n}$ for all $\left[ a_{ij}\right]
_{i,j=1}^{n}\in M_{n}\left( \mathcal{A}\right) $ such that $\left[ a_{ij}%
\right] _{i,j=1}^{n}=_{\lambda }\left[ 0\right] _{i,j=1}^{n}$ and for all $%
n, $ and so $\varphi $ is local completely positive.
\end{proof}

\begin{corollary}
\label{localpositive}Let $\varphi :\mathcal{A}\rightarrow $ $\mathcal{B}$ be
a linear map. Then $\varphi $ is local positive if and only if $\varphi $ is
continuous and positive.
\end{corollary}

\begin{remark}
\label{London} Let $\varphi :\mathcal{A}\rightarrow $ $\mathcal{B}$ be a
local completely positive map. If the locally $C^{\ast }$-algebra $\mathcal{A%
}$ is unital, then $\varphi $ is a strictly continuous completely positive\
map \cite[Remark 4.4]{J}. In particular, the structure theorem \cite[Theorem
4.6]{J} is valid for local completely positive maps on unital locally $%
C^{\ast }$-algebras.
\end{remark}

Let $\mathcal{CPCC}_{\text{loc}}(\mathcal{A},C^{\ast }(\mathcal{D}_{\mathcal{%
E}}))$ denote the set of all maps $\varphi :\mathcal{A}\rightarrow C^{\ast }(%
\mathcal{D}_{\mathcal{E}})$ which are local completely positive and local
completely contractive. If $\varphi \in \mathcal{CPCC}_{\text{loc}}(\mathcal{%
A},C^{\ast }(\mathcal{D}_{\mathcal{E}})),$ then $\varphi \left( b(\mathcal{A}%
)\right) \subseteq b(C^{\ast }(\mathcal{D}_{\mathcal{E}}))$. Moreover, there
is a completely contractive and completely positive map $\left. \varphi
\right\vert _{b(\mathcal{A})}:b(\mathcal{A})\rightarrow B(\mathcal{H})$ such
that $\left. \left. \varphi \right\vert _{b(\mathcal{A})}\left( a\right)
\right\vert _{\mathcal{D}_{\mathcal{E}}}=\varphi \left( a\right) \ $for all $%
a\in b(\mathcal{A}).$

The following result is a version of the Stinespring theorem for unbounded
local completely positive and local completely contractive maps.

\begin{theorem}
\label{s} \cite[Theorem 5.1]{D1} Let $\mathcal{A}$ be a unital locally $%
C^{\ast }$-algebra and $\varphi \in \mathcal{CPCC}_{\text{loc}}(\mathcal{A}%
,C^{\ast }(\mathcal{D}_{\mathcal{E}}))$. Then there exist a quantized domain 
$\{\mathcal{H}^{\varphi },\mathcal{E}^{\varphi },\mathcal{D}_{\mathcal{E}%
^{\varphi }}\}$, where $\mathcal{E}^{\varphi }=\{\mathcal{H}_{\lambda
}^{\varphi };\lambda \in \Lambda \}$ is an upward filtered family of closed
subspaces of $\mathcal{H}^{\varphi }$, a contraction $V_{\varphi }:\mathcal{H%
}\rightarrow \mathcal{H}^{\varphi }$ and a unital local contractive $\ast $%
-homomorphism $\pi _{\varphi }:\mathcal{A\rightarrow }C^{\ast }(\mathcal{D}_{%
\mathcal{E}^{\varphi }})$ such that

\begin{enumerate}
\item $V_{\varphi }\left( \mathcal{E}\right) \subseteq \mathcal{E}^{\varphi
};$

\item $\varphi \left( a\right) \subseteq V_{\varphi }^{\ast }\pi _{\varphi
}\left( a\right) V_{\varphi };$

\item $\mathcal{H}_{\lambda }^{\varphi }=\left[ \pi _{\varphi }\left( 
\mathcal{A}\right) V_{\varphi }\mathcal{H}_{\lambda }\right] $ for all $%
\lambda \in \Lambda .$

Moreover, if $\varphi \left( 1_{\mathcal{A}}\right) =$id$_{\mathcal{D}_{%
\mathcal{E}}}$, then $V_{\varphi }$ is an isometry.
\end{enumerate}
\end{theorem}

\ \ \ \ \ \ The triple $\left( \pi _{\varphi },V_{\varphi },\{\mathcal{H}%
^{\varphi },\mathcal{E}^{\varphi },\mathcal{D}_{\mathcal{E}^{\varphi
}}\}\right) $ constructed in Theorem \ref{s} is called a minimal Stinespring
dilation associated to $\varphi $. Moreover, the minimal Stinespring
dilation associated to $\varphi $ is unique up to unitary equivalence in the
following sense, if $(\ \pi _{\varphi },V_{\varphi },$\ $\{\mathcal{H}%
^{\varphi },\mathcal{E}^{\varphi },\mathcal{D}_{\mathcal{E}^{\varphi }}\})$
and$\left( \widetilde{\pi }_{\varphi },\widetilde{V}_{\varphi },\{\widetilde{%
\mathcal{H}}^{\varphi },\widetilde{\mathcal{E}}^{\varphi },\widetilde{%
\mathcal{D}}_{\widetilde{\mathcal{E}}^{\varphi }}\}\right) \ $are two
minimal Stinespring dilations associated to $\varphi $, then there is a
unitary operator $U_{\varphi }:\mathcal{H}^{\varphi }\rightarrow \widetilde{%
\mathcal{H}}^{\varphi }$ such that $U_{\varphi }V_{\varphi }=\widetilde{V}%
_{\varphi }$ and $U_{\varphi }\pi _{\varphi }\left( a\right) \subseteq 
\widetilde{\pi }_{\varphi }\left( a\right) U_{\varphi }$ for all $a\in 
\mathcal{A\ }$\cite[Theorem 3.4]{BGK}.

If $\left( \pi _{\varphi },V_{\varphi },\{\mathcal{H}^{\varphi },\mathcal{E}%
^{\varphi },\mathcal{D}_{\mathcal{E}^{\varphi }}\}\right) $ is a minimal
Stinespring dilation associated to $\varphi \in \mathcal{CPCC}_{\text{loc}}(%
\mathcal{A},C^{\ast }(\mathcal{D}_{\mathcal{E}}))$, then $\left( \left. \pi
_{\varphi }\right\vert _{b(\mathcal{A})},V_{\varphi },\mathcal{H}^{\varphi
}\right) $, where $\left. \left. \pi _{\varphi }\right\vert _{b(\mathcal{A}%
)}\left( a\right) \right\vert _{\mathcal{D}_{\widetilde{\mathcal{E}}%
^{\varphi }}}=\pi _{\varphi }\left( a\right) \ $for all $a\in b(\mathcal{A}%
), $ is a minimal Stinespring dilation associated with $\left. \varphi
\right\vert _{b(\mathcal{A})}.$

Let $\mathcal{E}=\{\mathcal{H}_{\iota };\iota \in \Upsilon \}$ and $\mathcal{%
F}=\{\mathcal{K}_{\gamma };\gamma \in \Gamma \}\ $be quantized domains in
Hilbert spaces $\mathcal{H}$ and $\mathcal{K}$, respectively. Then 
\begin{equation*}
\mathcal{E}\otimes \mathcal{F}=\left\{ \mathcal{H}_{\iota }\otimes \mathcal{K%
}_{\gamma };\left( \iota ,\gamma \right) \in \Upsilon \times \Gamma \right\}
\end{equation*}%
is a quantized domain in the Hilbert space $\mathcal{H}\otimes \mathcal{K}$,
with the union space $\mathcal{D}_{\mathcal{E}\otimes \mathcal{F}%
}=\tbigcup\limits_{\left( \iota ,\gamma \right) \in \Upsilon \times \Gamma }%
\mathcal{H}_{\iota }\otimes \mathcal{K}_{\gamma }.$\ If $\mathcal{E}=\{%
\mathcal{H}\}$, then $\mathcal{H}\otimes \mathcal{F}=\left\{ \mathcal{H}%
\otimes \mathcal{K}_{\gamma };\gamma \in \Gamma \right\} $ is a quantized
domain in the Hilbert space $\mathcal{H}\otimes \mathcal{K}$ with the union
space $\mathcal{D}_{\mathcal{H}\otimes \mathcal{F}}=\tbigcup\limits_{\gamma
\in \Gamma }\mathcal{H}\otimes \mathcal{K}_{\gamma }$.

The map $\Phi :C^{\ast }(\mathcal{D}_{\mathcal{E}})\otimes _{\text{alg}%
}C^{\ast }(\mathcal{D}_{\mathcal{F}})\rightarrow C^{\ast }(\mathcal{D}_{%
\mathcal{E}\otimes \mathcal{F}})$ given by 
\begin{equation*}
\Phi \left( T\otimes S\right) \left( \xi \otimes \eta \right) =T\xi \otimes
S\eta ,T\in C^{\ast }(\mathcal{D}_{\mathcal{E}}),S\in C^{\ast }(\mathcal{D}_{%
\mathcal{F}}),\xi \in \mathcal{D}_{\mathcal{E}},\eta \in \mathcal{D}_{%
\mathcal{F}}
\end{equation*}%
identifies $C^{\ast }(\mathcal{D}_{\mathcal{E}})\otimes _{\text{alg}}C^{\ast
}(\mathcal{D}_{\mathcal{F}})\ $with a $\ast $-subalgebra of $C^{\ast }(%
\mathcal{D}_{\mathcal{E}\otimes \mathcal{F}})$. The minimal tensor product
of the locally $C^{\ast }$-algebras $C^{\ast }(\mathcal{D}_{\mathcal{E}})$
and $C^{\ast }(\mathcal{D}_{\mathcal{F}})$ is the locally $C^{\ast }$%
-algebra $C^{\ast }(\mathcal{D}_{\mathcal{E}})\otimes C^{\ast }(\mathcal{D}_{%
\mathcal{F}})$ obtained by the completion of the $\ast $-subalgebra $C^{\ast
}(\mathcal{D}_{\mathcal{E}})\otimes _{\text{alg}}C^{\ast }(\mathcal{D}_{%
\mathcal{F}})$ in $C^{\ast }(\mathcal{D}_{\mathcal{E}\otimes \mathcal{F}})$
(see for example \cite{G}).

Let $\mathcal{A}$ and $\mathcal{B}$ be two locally $C^{\ast }$-algebras.
Recall that $\mathcal{A}$ and $\mathcal{B}$ can be identified with a locally 
$C^{\ast }$-subalgebra in $C^{\ast }(\mathcal{D}_{\mathcal{E}})$ and $%
C^{\ast }(\mathcal{D}_{\mathcal{F}})$, respectively for some quantized
domains $\mathcal{D}_{\mathcal{E}}$ and $\mathcal{D}_{\mathcal{F}}$. The
minimal or spatial tensor product of the locally $C^{\ast }$-algebras $%
\mathcal{A}$ and $\mathcal{B}$ is the locally $C^{\ast }$-algebra $\mathcal{A%
}\otimes \mathcal{B}$ obtained by the completion of the $\ast $-subalgebra $%
\mathcal{A}\otimes _{\text{alg}}\mathcal{B}$\ in $C^{\ast }(\mathcal{D}_{%
\mathcal{E}\otimes \mathcal{F}})$. In fact, $\mathcal{A}\otimes \mathcal{B}$
can be identified with the projective limit of the projective system of $%
C^{\ast }$-algebras $\{\mathcal{A}_{\lambda }\otimes \mathcal{B}_{\delta
};\pi _{\lambda _{1}\lambda _{2}}^{\mathcal{A}}\otimes \pi _{\delta
_{1}\delta _{2}}^{\mathcal{B}},\lambda _{1}\geq \lambda _{2},\delta _{1}\geq
\delta _{2}\}_{\left( \lambda ,\delta \right) \in \Lambda \times \Delta }$
(see \cite{Ph,F}).

If $\mathcal{A}$ is a $C^{\ast }$-algebra acting nondegenerately on a
Hilbert space $\mathcal{H}$ and $\mathcal{B}$ is a locally $C^{\ast }$%
-algebra that is identified with a locally $C^{\ast }$-subalgebra in $%
C^{\ast }(\mathcal{D}_{\mathcal{F}})$, then the minimal tensor product $%
\mathcal{A}\otimes \mathcal{B}$ of $\mathcal{A}$ and $\mathcal{B}$ is the
completion of the $\ast $-subalgebra $\mathcal{A}\otimes _{\text{alg}}%
\mathcal{B}$\ in $C^{\ast }(\mathcal{D}_{\mathcal{H}\otimes \mathcal{F}})$.

Let $\mathcal{A}$ and $\mathcal{B}$ be two locally $C^{\ast }$-algebras, $%
\pi _{1}:$ $\mathcal{A\rightarrow }C^{\ast }(\mathcal{D}_{\mathcal{E}})$ and 
$\pi _{2}:$ $\mathcal{B\rightarrow }C^{\ast }(\mathcal{D}_{\mathcal{F}})$ be
two local contractive $\ast $-morphisms. By the functorial property of the
minimal tensor product of $C^{\ast }$-algebras and taking into account the
above discussion we conclude that there is a unique local $\ast $-morphism $%
\pi _{1}\otimes \pi _{2}:\mathcal{A\otimes B\rightarrow }C^{\ast }(\mathcal{D%
}_{\mathcal{E}\otimes \mathcal{F}})$ such that 
\begin{equation*}
\left( \pi _{1}\otimes \pi _{2}\right) \left( a\otimes b\right) =\pi
_{1}\left( a\right) \otimes \pi _{2}\left( b\right)
\end{equation*}%
where $\left( \pi _{1}\left( a\right) \otimes \pi _{2}\left( b\right)
\right) \left( \xi \otimes \eta \right) =\pi _{1}\left( a\right) \xi \otimes
\pi _{2}\left( b\right) \eta $ for all $\xi \in \mathcal{D}_{\mathcal{E}%
},\eta \in \mathcal{D}_{\mathcal{F}}$.

\begin{proposition}
\label{3} Let $\mathcal{A}$ and $\mathcal{B}$ be two unital locally $C^{\ast
}$-algebras and let $\varphi \in \mathcal{CPCC}_{\text{loc}}(\mathcal{A}%
,C^{\ast }(\mathcal{D}_{\mathcal{E}}))$ and $\psi \in \mathcal{CPCC}_{\text{%
loc}}(\mathcal{B},C^{\ast }(\mathcal{D}_{\mathcal{F}}))$. Then there is a
unique map $\varphi \otimes \psi \in $ $\mathcal{CPCC}_{\text{loc}}(\mathcal{%
A\otimes B},C^{\ast }(\mathcal{D}_{\mathcal{E}\otimes \mathcal{F}}))$ such
that 
\begin{equation*}
\left( \varphi \otimes \psi \right) \left( a\otimes b\right) =\varphi \left(
a\right) \otimes \psi \left( b\right) \ \text{for all }a\in \mathcal{A\ }%
\text{and }b\in \mathcal{B}.
\end{equation*}%
Moreover, if $\left( \pi _{\varphi },V_{\varphi },\{\mathcal{H}^{\varphi },%
\mathcal{E}^{\varphi },\mathcal{D}_{\mathcal{E}^{\varphi }}\}\right) $ is a
minimal Stinespring dilation associated to $\varphi $ and $\left( \pi _{\psi
},V_{\psi },\{\mathcal{H}^{\psi },\mathcal{E}^{\psi },\mathcal{D}_{\mathcal{E%
}^{\psi }}\}\right) $ is a minimal Stinespring dilation associated to $\psi $%
, then $\left( \pi _{\varphi }\otimes \pi _{\psi },V_{\varphi }\otimes
V_{\psi },\{\mathcal{H}^{\varphi }\otimes \mathcal{H}^{\psi },\mathcal{E}%
^{\varphi }\otimes \mathcal{E}^{\psi },\mathcal{D}_{\mathcal{E}^{\varphi
}\otimes \mathcal{E}^{\psi }}\}\right) \ $is a minimal Stinespring dilation
associated to $\varphi \otimes \psi .$
\end{proposition}

\begin{proof}
Since $\pi _{\varphi }:\mathcal{A\rightarrow }C^{\ast }(\mathcal{D}_{%
\mathcal{E}^{\varphi }})$ and $\pi _{\psi }:\mathcal{B\rightarrow }C^{\ast }(%
\mathcal{D}_{\mathcal{E}^{\psi }})$ are local contractive $\ast $-morphisms,
there is a local contractive $\ast $-morphism $\pi _{\varphi }\otimes \pi
_{\psi }:\mathcal{A\otimes B\rightarrow }$ $C^{\ast }(\mathcal{D}_{\mathcal{E%
}^{\varphi }\otimes \mathcal{E}^{\psi }})$ such that 
\begin{equation*}
\left( \pi _{\varphi }\otimes \pi _{\psi }\right) \left( a\otimes b\right)
=\pi _{\varphi }\left( a\right) \otimes \pi _{\psi }\left( b\right) ,
\end{equation*}%
where $\left( \pi _{\varphi }\left( a\right) \otimes \pi _{\psi }\left(
b\right) \right) \left( \xi \otimes \eta \right) =\pi _{\varphi }\left(
a\right) \xi \otimes \pi _{\psi }\left( b\right) \eta $ for all $\xi \in 
\mathcal{D}_{\mathcal{E}^{\varphi }},\eta \in \mathcal{D}_{\mathcal{E}^{\psi
}}$. Since $V_{\varphi }:\mathcal{H}\rightarrow \mathcal{H}^{\varphi }$ and $%
V_{\psi }:\mathcal{K}\rightarrow \mathcal{H}^{\psi }$ are contractions, we
get a contraction $V_{\varphi }\otimes V_{\psi }:\mathcal{H}\otimes \mathcal{%
K}\rightarrow \mathcal{H}^{\varphi }\otimes \mathcal{H}^{\psi }$.$\ $Since $%
V_{\varphi }\left( \mathcal{E}\right) \subseteq \mathcal{E}^{\varphi }$ and $%
V_{\psi }\left( \mathcal{F}\right) \subseteq \mathcal{E}^{\psi },$ $\left(
V_{\varphi }\otimes V_{\psi }\right) \left( \mathcal{H}_{\iota }\otimes 
\mathcal{K}_{\mathcal{\gamma }}\right) \subseteq \mathcal{H}_{\iota
}^{\varphi }\otimes \mathcal{H}_{\mathcal{\gamma }}^{\varphi }\ $for all $%
\left( \iota ,\mathcal{\gamma }\right) \in \Upsilon \times \Gamma $, thus 
\begin{equation*}
\left( V_{\varphi }\otimes V_{\psi }\right) \left( \mathcal{E}\otimes 
\mathcal{F}\right) \subseteq \mathcal{E}^{\varphi }\otimes \mathcal{E}^{\psi
}.
\end{equation*}%
Therefore, we may consider the map $\varphi \otimes \psi :$ $\mathcal{%
A\otimes B\rightarrow }C^{\ast }(\mathcal{D}_{\mathcal{E}\otimes \mathcal{F}%
})$ given by%
\begin{equation*}
\left( \varphi \otimes \psi \right) \left( c\right) =\left. \left(
V_{\varphi }\otimes V_{\psi }\right) ^{\ast }\left( \pi _{\varphi }\otimes
\pi _{\psi }\right) \left( c\right) \left( V_{\varphi }\otimes V_{\psi
}\right) \right\vert _{\mathcal{D}_{\mathcal{E}^{\varphi }\otimes \mathcal{E}%
^{\psi }}}\text{.}
\end{equation*}%
Since $\pi _{\varphi }\otimes \pi _{\psi }$ is a local contractive $\ast $%
-morphism and $V_{\varphi }\otimes V_{\psi }$ is a contraction, $\varphi
\otimes \psi $ is a local completely positive and local completely
contractive map from $\mathcal{A\otimes B}$ to $C^{\ast }(\mathcal{D}_{%
\mathcal{E}\otimes \mathcal{F}})$. Clearly, $\left( \varphi \otimes \psi
\right) \left( a\otimes b\right) =\varphi \left( a\right) \otimes \psi
\left( b\right) \ $for all $a\in \mathcal{A\ }$and $b\in \mathcal{B}$.

To show that $\left( \pi _{\varphi }\otimes \pi _{\psi },V_{\varphi }\otimes
V_{\psi },\{\mathcal{H}^{\varphi }\otimes \mathcal{H}^{\psi },\mathcal{E}%
^{\varphi }\otimes \mathcal{E}^{\psi },\mathcal{D}_{\mathcal{E}^{\varphi
}\otimes \mathcal{E}^{\psi }}\}\right) \ $is a minimal Stinespring dilation
associated with $\varphi \otimes \psi $ it remains to show that $\mathcal{H}%
_{\iota }^{\varphi }\otimes \mathcal{H}_{\gamma }^{\psi }=\left[ \left( \pi
_{\varphi }\otimes \pi _{\psi }\right) \left( \mathcal{A\otimes B}\right)
\left( V_{\varphi }\otimes V_{\psi }\right) \left( \mathcal{H}_{\iota
}\otimes \mathcal{K}_{\gamma }\right) \right] $ for all $\left( \iota ,%
\mathcal{\gamma }\right) \in \Upsilon \times \Gamma $.

Let $\left( \iota ,\mathcal{\gamma }\right) \in \Upsilon \times \Gamma $.
Then 
\begin{eqnarray*}
\left[ \left( \pi _{\varphi }\otimes \pi _{\psi }\right) \left( \mathcal{%
A\otimes B}\right) \left( V_{\varphi }\otimes V_{\psi }\right) \left( 
\mathcal{H}_{\iota }\otimes \mathcal{K}_{\gamma }\right) \right] &=&\left[
\pi _{\varphi }\left( \mathcal{A}\right) V_{\varphi }\left( \mathcal{H}%
_{\iota }\right) \otimes \pi _{\psi }\left( \mathcal{B}\right) V_{\psi
}\left( \mathcal{K}_{\gamma }\right) \right] \\
&=&\mathcal{H}_{\iota }^{\varphi }\otimes \mathcal{H}_{\gamma }^{\psi },
\end{eqnarray*}%
as required.

To show the uniqueness of the map $\varphi \otimes \psi $, let $\phi \in 
\mathcal{CPCC}_{\text{loc}}(\mathcal{A\otimes B},C^{\ast }(\mathcal{D}_{%
\mathcal{E}\otimes \mathcal{F}}))$ such that $\phi \left( a\otimes b\right)
=\varphi \left( a\right) \otimes \psi \left( b\right) \ $for all $a\in 
\mathcal{A\ }$and $b\in \mathcal{B}$, $\phi \neq \varphi \otimes \psi $.
Since $\varphi \otimes \psi $ and $\phi $ are continuous and $\mathcal{%
A\otimes }_{\text{alg}}\mathcal{B}$ is dense, from $\phi \left( a\otimes
b\right) =\varphi \left( a\right) \otimes \psi \left( b\right) =\left(
\varphi \otimes \psi \right) \left( a\otimes b\right) $ for all $a\in 
\mathcal{A\ }$and $b\in \mathcal{B}$ it follows that $\phi =\varphi \otimes
\psi $, a contradiction, and the uniqueness of the map $\varphi \otimes \psi 
$ is proved.
\end{proof}

\section{Main results}

\ Suppose that $\{\mathcal{H};\mathcal{E};\mathcal{D}_{\mathcal{E}}\}$ is a
Fr\'{e}chet quantized domain in a Hilbert space $\mathcal{H}\ $(that is, $%
\mathcal{E=\{H}_{n}\mathcal{\}}_{n\in \mathbb{N}}$). For each $n\in \mathbb{N%
}^{\ast },$ let $\mathcal{H}_{n}^{c}$ be the orthogonal complement of the
closed subspace\ $\mathcal{H}_{n-1}$ in $\mathcal{H}_{n}$, and put $\mathcal{%
H}_{0}^{c}=\mathcal{H}_{0}$. Then $\mathcal{H}_{n}=\tbigoplus\limits_{k\leq
n}\mathcal{H}_{k}^{c}$. For each $n\in \mathbb{N}$ and for each $\xi \in 
\mathcal{H}_{k}^{c},k\leq n$, the rank one operator $\theta _{\xi ,\xi }:%
\mathcal{D}_{\mathcal{E}}\rightarrow \mathcal{D}_{\mathcal{E}},$ $\theta
_{\xi ,\xi }(\eta )=\xi \left\langle \xi ,\eta \right\rangle \ $is an
element in $C^{\ast }(\mathcal{D}_{\mathcal{E}})$. The closed two sided $%
\ast $-ideal of $C^{\ast }(\mathcal{D}_{\mathcal{E}})$ generated by these
rank one operators is denoted by $K(\mathcal{D}_{\mathcal{E}})$ and it is
called the locally $C^{\ast }$-algebra of all compact operators on $\mathcal{%
D}_{\mathcal{E}}$. Clearly, for each $n\in \mathbb{N},K(\mathcal{D}_{%
\mathcal{E}})_{n}$ $\subseteq K\left( \mathcal{H}_{n}\right) $, which it is
a closed two sided $\ast $-ideal of $C^{\ast }(\mathcal{D}_{\mathcal{E}%
})_{n}.$

Let $\mathcal{A}$ be a unital locally $C^{\ast }$-algebra, $\{\mathcal{H};%
\mathcal{E};\mathcal{D}_{\mathcal{E}}\}$ be a Fr\'{e}chet quantized domain
in a Hilbert space $\mathcal{H}$ with $\mathcal{E=\{H}_{n}\mathcal{\}}_{n\in 
\mathbb{N}}$, $a\in \mathcal{A}$ and $c\in \mathcal{A\otimes }C^{\ast }(%
\mathcal{D}_{\mathcal{E}})$. If $\xi \in \mathcal{H}_{n}^{c}$ and $a\otimes
\theta _{\xi ,\xi }\geq _{\left( \lambda ,n\right) }c$ $\geq _{\left(
\lambda ,n\right) }0$, then there is $b\in \mathcal{A}$ such that $%
c=_{\left( \lambda ,n\right) }b\otimes \theta _{\xi ,\xi }\ $. Indeed, from $%
a\otimes \theta _{\xi ,\xi }\geq _{\left( \lambda ,n\right) }c$ $\geq
_{\left( \lambda ,n\right) }0$ we deduce that 
\begin{equation*}
\pi _{\lambda }^{\mathcal{A}}\left( a\right) \otimes \left. \theta _{\xi
,\xi }\right\vert _{\mathcal{H}_{n}}=\pi _{\left( \lambda ,n\right) }^{%
\mathcal{A\otimes }C^{\ast }(\mathcal{D}_{\mathcal{E}})}\left( a\otimes
\theta _{\xi ,\xi }\right) \geq \pi _{\left( \lambda ,n\right) }^{\mathcal{%
A\otimes }C^{\ast }(\mathcal{D}_{\mathcal{E}})}\left( c\right) \geq 0
\end{equation*}%
and by \cite[Lemma 2.1]{BO}, there is $b_{\lambda }\in \mathcal{A}_{\lambda
} $ such that $\pi _{\left( \lambda ,n\right) }^{\mathcal{A\otimes }C^{\ast
}(\mathcal{D}_{\mathcal{E}})}\left( c\right) =b_{\lambda }\otimes \left.
\theta _{\xi ,\xi }\right\vert _{\mathcal{H}_{n}}.$ Therefore, there is $%
b\in \mathcal{A}$ such that 
\begin{equation*}
\pi _{\left( \lambda ,n\right) }^{\mathcal{A\otimes }C^{\ast }(\mathcal{D}_{%
\mathcal{E}})}\left( c\right) =b_{\lambda }\otimes \left. \theta _{\xi ,\xi
}\right\vert _{\mathcal{H}_{n}}=\pi _{\lambda }^{\mathcal{A}}\left( b\right)
\otimes \left. \theta _{\xi ,\xi }\right\vert _{\mathcal{H}_{n}}=\pi
_{\left( \lambda ,n\right) }^{\mathcal{A\otimes }C^{\ast }(\mathcal{D}_{%
\mathcal{E}})}\left( b\otimes \theta _{\xi ,\xi }\right) .
\end{equation*}

\begin{theorem}
\label{6} Let $\mathcal{A}$ and $\mathcal{B}$\ be two unital locally $%
C^{\ast }$-algebras, $\{\mathcal{H};\mathcal{E};\mathcal{D}_{\mathcal{E}}\}$
be a Fr\'{e}chet quantized domain in a Hilbert space $\mathcal{H}$ with $%
\mathcal{E=\{H}_{n}\mathcal{\}}_{n\in \mathbb{N}},$ $\varphi :$ $\mathcal{A}$
$\rightarrow $ $\mathcal{B\ }$and $\phi :\mathcal{A\otimes }C^{\ast }(%
\mathcal{D}_{\mathcal{E}})\rightarrow \mathcal{B\otimes }C^{\ast }(\mathcal{D%
}_{\mathcal{E}})$ be two linear maps. Then the following conditions are
equivalent:

\begin{enumerate}
\item $\phi $ and $\varphi \otimes $id$_{C^{\ast }(\mathcal{D}_{\mathcal{E}%
})}-\phi $ are local positive.

\item There exists a local positive map $\psi :\mathcal{A\rightarrow B\ }$%
such that $\varphi -\psi $ is local positive and $\phi =\psi \otimes $ id$%
_{C^{\ast }(\mathcal{D}_{\mathcal{E}})}$.
\end{enumerate}
\end{theorem}

\begin{proof}
$(1)\Rightarrow (2)$ Since the maps $\phi $ and $\varphi \otimes $id$%
_{C^{\ast }(\mathcal{D}_{\mathcal{E}})}-\phi $ are local positive, the map $%
\varphi \otimes $id$_{C^{\ast }(\mathcal{D}_{\mathcal{E}})}$ is local
positive. Then, $\varphi $ is local positive, and by Corollary \ref%
{localpositive} 
\begin{equation*}
0\leq \phi \left( a\otimes T\right) \leq \varphi \left( a\right) \otimes T
\end{equation*}%
for all $a\in \mathcal{A}\ $with $a\geq 0$ and for all $T\in C^{\ast }(%
\mathcal{D}_{\mathcal{E}})$ with $T\geq 0.$ Therefore, for each $\left(
\delta ,n\right) \in \Delta \times \mathbb{N},$ there are $\lambda _{0}\in
\Lambda $ and $C_{0}>0$ such that%
\begin{equation*}
\left\Vert \phi \left( a\otimes T\right) \right\Vert _{\left( \delta
,n\right) }\leq \left\Vert \varphi \left( a\right) \otimes T\right\Vert
_{\left( \delta ,n\right) }=\left\Vert \varphi \left( a\right) \right\Vert
_{\delta }\left\Vert T\right\Vert _{n}\leq C_{0}p_{\lambda _{0}}\left(
a\right) \left\Vert T\right\Vert _{n}
\end{equation*}%
for all $a\in \mathcal{A}\ $with $a\geq 0$ and for all $T\in C^{\ast }(%
\mathcal{D}_{\mathcal{E}})$ with $T\geq 0$. Therefore, there exist maps $%
\phi _{\left( \delta ,n\right) }:$ $\mathcal{A}_{\lambda _{0}}\mathcal{%
\otimes }C^{\ast }(\mathcal{D}_{\mathcal{E}})_{n}\rightarrow \mathcal{B}%
_{\delta }\mathcal{\otimes }C^{\ast }(\mathcal{D}_{\mathcal{E}})_{n}$ and $%
\varphi _{\delta }:\mathcal{A}_{\lambda _{0}}\rightarrow \mathcal{B}_{\delta
}$ such that 
\begin{equation*}
\pi _{\left( \delta ,n\right) }^{\mathcal{B\otimes }C^{\ast }(\mathcal{D}_{%
\mathcal{E}})}\circ \phi =\phi _{\left( \delta ,n\right) }\circ \pi _{\left(
\lambda _{0},n\right) }^{\mathcal{A\otimes }C^{\ast }(\mathcal{D}_{\mathcal{E%
}})}\text{ and }\pi _{\delta }^{\mathcal{B}}\circ \varphi =\varphi _{\delta
}\circ \pi _{\lambda _{0}}^{\mathcal{A}}.
\end{equation*}%
Moreover, $\phi _{\left( \delta ,n\right) }$ and $\varphi _{\delta }\otimes $%
id$_{C^{\ast }(\mathcal{D}_{\mathcal{E}})_{n}}-\phi _{\left( \delta
,n\right) }$ are positive. Thus, by the proof of \cite[Theorem 2]{BO}, there
is a positive map $\widetilde{\psi }_{\delta }:\mathcal{A}_{\lambda
_{0}}\rightarrow \mathcal{B}_{\delta }$ such that 
\begin{equation*}
\phi _{\left( \delta ,n\right) }\left( \pi _{\lambda _{0}}^{\mathcal{A}%
}\left( a\right) \otimes \left. \theta _{\xi ,\xi }\right\vert _{\mathcal{H}%
_{n}}\right) =\widetilde{\psi }_{\delta }\left( \pi _{\lambda _{0}}^{%
\mathcal{A}}\left( a\right) \right) \otimes \left. \theta _{\xi ,\xi
}\right\vert _{\mathcal{H}_{n}}
\end{equation*}%
for all $a\in \mathcal{A}$ such that $a\geq _{\lambda _{0}}0$ and $\xi \in 
\mathcal{H}_{k}^{c},k\leq n$, and $\varphi _{\delta }-\widetilde{\psi }%
_{\delta }$ is positive. Since $K(\mathcal{D}_{\mathcal{E}})_{n}\ $is
generated by rank one operators $\theta _{\xi ,\xi },\xi \in \mathcal{H}%
_{k}^{c},$ $k\leq n,$ 
\begin{equation*}
\phi _{\left( \delta ,n\right) }\left( \pi _{\lambda _{0}}^{\mathcal{A}%
}\left( a\right) \otimes \left. T\right\vert _{\mathcal{H}_{n}}\right) =%
\widetilde{\psi }_{\delta }\left( \pi _{\lambda _{0}}^{\mathcal{A}}\left(
a\right) \right) \otimes \left. T\right\vert _{\mathcal{H}_{n}}
\end{equation*}%
for all $a\in \mathcal{A}$ and for all $T\in K(\mathcal{D}_{\mathcal{E}}).$

Suppose that $\mathcal{B}_{\delta }$ acts nondegenerately on a Hilbert space 
$\mathcal{K}$. Let $\{u_{i}\}_{i\in I}$ be an approximate unit for $K(%
\mathcal{D}_{\mathcal{E}})_{n}$, $\eta \in \mathcal{K}$ and $\xi \in 
\mathcal{H}_{k}^{c},k\leq n,\left\Vert \xi \right\Vert \neq 0.$ Since $%
\mathcal{H}_{n}=\tbigoplus\limits_{k\leq n}\mathcal{H}_{k}^{c}$ and $K(%
\mathcal{D}_{\mathcal{E}})_{n}$ is a closed two sided $\ast $-ideal of $%
C^{\ast }(\mathcal{D}_{\mathcal{E}})_{n}$, we have 
\begin{eqnarray*}
&&\left\langle \phi _{\left( \delta ,n\right) }\left( \pi _{\lambda _{0}}^{%
\mathcal{A}}\left( a\right) \otimes \left. T\right\vert _{\mathcal{H}%
_{n}}\right) \left( \eta \otimes \xi \right) ,\left( \eta \otimes \xi
\right) \right\rangle \\
&=&\lim\limits_{i}\left\langle \phi _{\left( \delta ,n\right) }\left( \pi
_{\lambda _{0}}^{\mathcal{A}}\left( a\right) \otimes u_{i}\left.
T\right\vert _{\mathcal{H}_{n}}u_{i}\right) \left( \eta \otimes \xi \right)
,\eta \otimes \xi \right\rangle \\
&=&\lim\limits_{i}\left\langle \left( \widetilde{\psi }_{\delta }\left( \pi
_{\lambda _{0}}^{\mathcal{A}}\left( a\right) \right) \otimes u_{i}\left.
T\right\vert _{\mathcal{H}_{n}}u_{i}\right) \left( \eta \otimes \xi \right)
,\eta \otimes \xi \right\rangle \\
&=&\lim\limits_{i}\left( \left\langle \widetilde{\psi }_{\delta }\left( \pi
_{\lambda _{0}}^{\mathcal{A}}\left( a\right) \right) \eta ,\eta
\right\rangle \otimes \left\langle u_{i}\left. T\right\vert _{\mathcal{H}%
_{n}}u_{i}\left. \theta _{\frac{1}{\left\Vert \xi \right\Vert }\xi ,\frac{1}{%
\left\Vert \xi \right\Vert }\xi }\right\vert _{\mathcal{H}_{n}}\left( \xi
\right) ,\xi \right\rangle \right) \\
&=&\left\langle \widetilde{\psi }_{\delta }\left( \pi _{\lambda _{0}}^{%
\mathcal{A}}\left( a\right) \right) \eta ,\eta \right\rangle \otimes
\left\langle \left. T\right\vert _{\mathcal{H}_{n}}\left. \theta _{\frac{1}{%
\left\Vert \xi \right\Vert }\xi ,\frac{1}{\left\Vert \xi \right\Vert }\xi
}\right\vert _{\mathcal{H}_{n}}\left( \xi \right) ,\xi \right\rangle \\
&=&\left\langle \widetilde{\psi }_{\delta }\left( \pi _{\lambda _{0}}^{%
\mathcal{A}}\left( a\right) \right) \eta ,\eta \right\rangle \otimes
\left\langle \left. T\right\vert _{\mathcal{H}_{n}}\left( \xi \right) ,\xi
\right\rangle \\
&=&\left\langle \left( \widetilde{\psi }_{\delta }\left( \pi _{\lambda
_{0}}^{\mathcal{A}}\left( a\right) \right) \otimes \left. T\right\vert _{%
\mathcal{H}_{n}}\right) \left( \eta \otimes \xi \right) ,\eta \otimes \xi
\right\rangle
\end{eqnarray*}%
for all $a\in \mathcal{A}$ and for all $T\in C^{\ast }(\mathcal{D}_{\mathcal{%
E}})$. Therefore, 
\begin{equation*}
\phi _{\left( \delta ,n\right) }\left( \pi _{\lambda _{0}}^{\mathcal{A}%
}\left( a\right) \otimes \left. T\right\vert _{\mathcal{H}_{n}}\right) =%
\widetilde{\psi }_{\delta }\left( \pi _{\lambda _{0}}^{\mathcal{A}}\left(
a\right) \right) \otimes \left. T\right\vert _{\mathcal{H}_{n}}
\end{equation*}%
for all $a\in \mathcal{A}$ and for all $T\in C^{\ast }(\mathcal{D}_{\mathcal{%
E}}).$

Let $\psi _{\delta }:\mathcal{A}\rightarrow \mathcal{B}_{\delta },\psi
_{\delta }=\widetilde{\psi }_{\delta }\circ \pi _{\lambda _{0}}^{\mathcal{A}%
} $. Clearly, $\psi _{\delta }$ is a local positive map, and 
\begin{equation*}
\pi _{\left( \delta ,n\right) }^{\mathcal{B\otimes }C^{\ast }(\mathcal{D}_{%
\mathcal{E}})}\left( \phi \left( a\otimes T\right) \right) =\psi _{\delta
}\left( a\right) \otimes \left. T\right\vert _{\mathcal{H}_{n}}
\end{equation*}%
for all $a\in \mathcal{A}$ and $T\in $ $C^{\ast }(\mathcal{D}_{\mathcal{E}})$%
.

Let $\delta _{1},\delta _{2}\in \Delta $ with $\delta _{1}\geq \delta _{2},$ 
$n\in \mathbb{N}$ and $a\in \mathcal{A}$. Since 
\begin{eqnarray*}
\pi _{\delta _{1}\delta _{2}}^{\mathcal{B}}\left( \psi _{\delta _{1}}\left(
a\right) \right) \otimes \left. T\right\vert _{\mathcal{H}_{n}} &=&\pi
_{\left( \delta _{1},n\right) \left( \delta _{2},n\right) }^{\mathcal{%
B\otimes }C^{\ast }(\mathcal{D}_{\mathcal{E}})}\left( \psi _{\delta
_{1}}\left( a\right) \otimes \left. T\right\vert _{\mathcal{H}_{n}}\right) \\
&=&\pi _{\left( \delta _{1},n\right) \left( \delta _{2},n\right) }^{\mathcal{%
B\otimes }C^{\ast }(\mathcal{D}_{\mathcal{E}})}\left( \pi _{\left( \delta
_{1},n\right) }^{\mathcal{B\otimes }C^{\ast }(\mathcal{D}_{\mathcal{E}%
})}\left( \phi \left( a\otimes T\right) \right) \right) \\
&=&\pi _{\left( \delta _{2},n\right) }^{\mathcal{B\otimes }C^{\ast }(%
\mathcal{D}_{\mathcal{E}})}\left( \phi \left( a\otimes T\right) \right)
=\psi _{\delta _{2}}\left( a\right) \otimes \left. T\right\vert _{\mathcal{H}%
_{n}}
\end{eqnarray*}%
for all $T\in $ $C^{\ast }(\mathcal{D}_{\mathcal{E}})$, it follows that $\pi
_{\delta _{1}\delta _{2}}^{\mathcal{B}}\left( \psi _{\delta _{1}}\left(
a\right) \right) =\psi _{\delta _{2}}\left( a\right) $. Therefore, there is
a linear map $\psi :\mathcal{A\rightarrow B}$ such that 
\begin{equation*}
\psi \left( a\right) =\left( \psi _{\delta }\left( a\right) \right) _{\delta
\in \Delta }\text{.}
\end{equation*}%
Moreover, since for each $\delta \in \Delta $, $\psi _{\delta }$ is local
positive, $\psi $ is local positive and 
\begin{eqnarray*}
\pi _{\left( \delta ,n\right) }^{\mathcal{B\otimes }C^{\ast }(\mathcal{D}_{%
\mathcal{E}})}\left( \left( \psi \otimes \text{id}_{C^{\ast }(\mathcal{D}_{%
\mathcal{E}})}\right) \left( a\otimes T\right) \right) &=&\psi _{\delta
}\left( a\right) \otimes \left. T\right\vert _{\mathcal{H}_{n}} \\
&=&\pi _{\left( \delta ,n\right) }^{\mathcal{B\otimes }C^{\ast }(\mathcal{D}%
_{\mathcal{E}})}\left( \phi \left( a\otimes T\right) \right)
\end{eqnarray*}%
for all $a\in \mathcal{A}$ and $T\in $ $C^{\ast }(\mathcal{D}_{\mathcal{E}})$%
, and for all $\left( \delta ,n\right) \in \Delta \times \mathbb{N}$.
Therefore, 
\begin{equation*}
\phi =\psi \otimes \text{id}_{C^{\ast }(\mathcal{D}_{\mathcal{E}})}\text{.}
\end{equation*}

To show that $\varphi -\psi $ is local positive, let $\delta \in \Delta $.
We seen that there exist $\lambda _{0}\in \Lambda $ such that $\varphi
_{\delta }-\widetilde{\psi }_{\delta }$ is positive, where $\pi _{\delta }^{%
\mathcal{B}}\circ \varphi =\varphi _{\delta }\circ \pi _{\lambda _{0}}^{%
\mathcal{A}}$ and $\pi _{\delta }^{\mathcal{B}}\circ \psi =\psi _{\delta }=%
\widetilde{\psi }_{\delta }\circ \pi _{\lambda _{0}}^{\mathcal{A}}$.
Clearly, $\left( \varphi -\psi \right) \left( a\right) \geq _{\delta }0$
whenever $a\geq _{\lambda _{0}}0$ and $\left( \varphi -\psi \right) \left(
a\right) =_{\delta }0$ whenever $a=_{\lambda _{0}}0$. Therefore, $\varphi
-\psi $ is local positive.

$\left( 2\right) \Rightarrow \left( 1\right) $ By Proposition \ref{3}, $\phi 
$ is local positive, and since 
\begin{equation*}
\varphi \otimes \text{id}_{C^{\ast }(\mathcal{D}_{\mathcal{E}})}-\phi
=\left( \varphi -\psi \right) \otimes \text{id}_{C^{\ast }(\mathcal{D}_{%
\mathcal{E}})}
\end{equation*}%
and $\varphi -\psi $ are local positive, $\varphi \otimes $id$_{C^{\ast }(%
\mathcal{D}_{\mathcal{E}})}-\phi $ is local positive.
\end{proof}

\begin{corollary}
\label{positive} Let $\mathcal{A}$ and $\mathcal{B}$\ be two unital locally $%
C^{\ast }$-algebras, $\mathcal{H}$ be a Hilbert space $\mathcal{H},$ $%
\varphi :$ $\mathcal{A}$ $\rightarrow $ $\mathcal{B\ }$and $\phi :\mathcal{%
A\otimes }B(\mathcal{H})\rightarrow \mathcal{B\otimes }B(\mathcal{H})$ be
two linear maps. Then $\phi $ and $\varphi \otimes $id$_{B(\mathcal{H}%
))}-\phi $ are local positive if and only if there is a local positive map $%
\psi :\mathcal{A\rightarrow B\ }$such that $\varphi -\psi $ is local
positive and $\phi =\psi \otimes $ id$_{B(\mathcal{H})}$.
\end{corollary}

As an application of Theorem \ref{6}, we show that given a local positive
map $\varphi :$ $\mathcal{A}\rightarrow $ $\mathcal{B}$, the local positive
map $\varphi \otimes $id$_{M_{n}\left( \mathbb{C}\right) }$ is local
decomposable for some $n\geq 2$ if and only if $\varphi $ is a local\ $CP$%
-map.

\begin{definition}
\cite{MJ3} A linear map $\varphi :\mathcal{A}\rightarrow \mathcal{B}$ is
called \textit{local }$\mathit{n}$\textit{-copositive if }the map $\varphi
\otimes t:\mathcal{A\otimes }M_{n}\left( \mathbb{C}\right) \rightarrow 
\mathcal{B\otimes }M_{n}\left( \mathbb{C}\right) $ defined by 
\begin{equation*}
\left( \varphi \otimes t\right) \left( \left[ a_{ij}\right]
_{i,j=1}^{n}\right) =\left[ \varphi \left( a_{ji}\right) \right] _{i,j=1}^{n}
\end{equation*}%
is local positive, where $t$ denotes the transpose map on $M_{n}\left( 
\mathbb{C}\right) $.

We say that $\varphi $ is local completely copositive if for each $\delta
\in \Delta $, there exists $\lambda \in \Lambda $ such that $\left( \varphi
\otimes t\right) \left( \left[ a_{ij}\right] _{i,j=1}^{n}\right) \geq
_{\delta }0\ $whenever $\left[ a_{ij}\right] _{i,j=1}^{n}\geq _{\lambda }0$
and $\left( \varphi \otimes t\right) \left( \left[ a_{ij}\right]
_{i,j=1}^{n}\right) =_{\delta }0\ \ $whenever $\left[ a_{ij}\right]
_{i,j=1}^{n}=_{\lambda }0$,$\ $for all $n\in \mathbb{N}.$
\end{definition}

\begin{remark}
\cite{MJ3}\label{copositive}Let $\mathcal{A}$ and $\mathcal{B}$\ be two
locally $C^{\ast }$-algebras and $\varphi :\mathcal{A}\rightarrow \mathcal{B}
$ be a \textit{local }$\mathit{n}$\textit{-copositive} map. Then:

\begin{enumerate}
\item $\varphi $ is \textit{\ }local positive and so it is continuous and
positive.

\item \textit{for each }$\delta \in \Delta $, there exist $\lambda \in
\Lambda $ and an $n$-copositive map $\varphi _{\delta }:\mathcal{A}_{\lambda
}\rightarrow \mathcal{B}_{\delta }$ such that $\pi _{\delta }^{\mathcal{B}%
}\circ \varphi =\varphi _{\delta }\circ \pi _{\lambda }^{\mathcal{A}}.$
\end{enumerate}
\end{remark}

\begin{definition}
\cite{MJ3} A linear map $\varphi :\mathcal{A}\rightarrow \mathcal{B}$ is
local decomposable if it is sum of a local completely positive map and a 
\textit{local completely copositive }map.
\end{definition}

If $\varphi :\mathcal{A}\rightarrow \mathcal{B}$ is local decomposable, then
for each\textit{\ }$\delta \in \Delta $, there exist $\lambda \in \Lambda $
and a decomposable positive map $\varphi _{\delta }:\mathcal{A}_{\lambda
}\rightarrow \mathcal{B}_{\delta }$ such that $\pi _{\delta }^{\mathcal{B}%
}\circ \varphi =\varphi _{\delta }\circ \pi _{\lambda }^{\mathcal{A}}.$

\begin{theorem}
\label{7}Let $\mathcal{A}$ and $\mathcal{B}$\ be two unital locally $C^{\ast
}$-algebras and $\varphi :\mathcal{A\rightarrow B}$ be a linear map. If for
some $n\geq 2$, $\varphi \otimes $id$_{M_{n}\left( \mathbb{C}\right) }:%
\mathcal{A\otimes }M_{n}\left( \mathbb{C}\right) \rightarrow \mathcal{%
B\otimes }M_{n}\left( \mathbb{C}\right) $ is local decomposable, then $%
\varphi $ is local completely positive.
\end{theorem}

\begin{proof}
We adapt the proof of \cite[Theorem 3.1 ]{BO}. If $\varphi \otimes $id$%
_{M_{n}\left( \mathbb{C}\right) }$ is decomposable, there are a local
completely positive map $\phi $ and a local completely copositive map $\psi $
such that $\varphi \otimes $id$_{M_{n}\left( \mathbb{C}\right) }=\phi +\psi $%
.$\ $By Corollary \ref{positive}, there exist two local positive maps $\phi
_{1}:$ $\mathcal{A\rightarrow B}$ and $\psi _{1}:$ $\mathcal{A\rightarrow B\ 
}$such that $\phi =\phi _{1}\otimes $id$_{M_{n}\left( \mathbb{C}\right) }$
and $\psi =\psi _{1}\otimes $id$_{M_{n}\left( \mathbb{C}\right) }$.

Since $\phi =\phi _{1}\otimes $id$_{M_{n}\left( \mathbb{C}\right) }$ and $%
\phi $ is local completely positive, $\phi _{1}$ is local completely
positive.

Since $\psi $ and $\psi _{1}$ are continuous (as $\psi $ is local completely
copositive and $\psi _{1}$ is local positive) and $\psi =\psi _{1}\otimes $id%
$_{M_{n}\left( \mathbb{C}\right) }$,\ for each $\delta \in \Delta ,$ there
exist $\lambda \in \Lambda $, and completely positive maps $\psi _{\delta }:%
\mathcal{A}_{\lambda }\rightarrow \mathcal{B}_{\delta }$ and $\psi _{1\delta
}:$ $\mathcal{A}_{\lambda }\rightarrow \mathcal{B}_{\delta }$\ such that $%
\pi _{\delta }^{\mathcal{B}}\circ \psi =\psi _{\delta }\circ \pi _{\lambda
}^{\mathcal{A}},$ $\pi _{\delta }^{\mathcal{B}}\circ \psi _{1}=\psi
_{1\delta }\circ \pi _{\lambda }^{\mathcal{A}}\ $and $\psi _{\delta }=\psi
_{1\delta }\otimes $id$_{M_{n}\left( \mathbb{C}\right) }$. Then, since $%
n\geq 2$, by \cite[Lemma 3.2]{BO}, for each $\delta \in \Delta ,$ $\psi
_{\delta }=0$. Consequently, $\psi =0$, and so $\varphi \otimes $id$%
_{M_{n}\left( \mathbb{C}\right) }=\phi =\phi _{1}\otimes $id$_{M_{n}\left( 
\mathbb{C}\right) }$, whence $\varphi =\phi _{1}.$
\end{proof}

Let $\{\mathcal{H};\mathcal{E};\mathcal{D}_{\mathcal{E}}\}$ be a quantized
domain in the Hilbert space $\mathcal{H\ }$with $\mathcal{E=\{H}_{\iota }%
\mathcal{\}}_{\iota \in \Upsilon }$.

For a local contractive $\ast $-morphism $\pi :\mathcal{A\rightarrow }%
C^{\ast }(\mathcal{D}_{\mathcal{E}})$ 
\begin{equation*}
\pi \left( \mathcal{A}\right) ^{^{\prime }}=\{T\in B\left( \mathcal{H}%
\right) ;T\pi \left( a\right) \subseteq \pi \left( a\right) T\text{ for all }%
a\in \mathcal{A}\}.
\end{equation*}

\begin{remark}
\begin{enumerate}
\item $\pi \left( \mathcal{A}\right) ^{^{\prime }}$is a\textbf{\ }von
Neumann algebra. Indeed,%
\begin{eqnarray*}
\pi \left( \mathcal{A}\right) ^{^{\prime }} &=&\{T\in B\left( \mathcal{H}%
\right) ;T\pi \left( a\right) \subseteq \pi \left( a\right) T\text{ for all }%
a\in \mathcal{A}\} \\
&=&\{T\in B\left( \mathcal{H}\right) ;T\pi \left( a\right) \subseteq \pi
\left( a\right) T\text{ for all }a\in b(\mathcal{A)}\} \\
&=&\{T\in B\left( \mathcal{H}\right) ;T\left. \pi \right\vert _{b(\mathcal{A)%
}}\left( a\right) =\left. \pi \right\vert _{b(\mathcal{A)}}\left( a\right) T%
\text{ for all }a\in b(\mathcal{A)}\} \\
&=&\left. \pi \right\vert _{b(\mathcal{A)}}\left( b(\mathcal{A})\right)
^{\prime }
\end{eqnarray*}%
where $\left. \pi \right\vert _{b(\mathcal{A)}}:b(\mathcal{A)\rightarrow }%
B\left( \mathcal{H}\right) $ is the $\ast $-representation of the $C^{\ast }$%
-algebra $b(\mathcal{A)}$ of all bounded elements of $\mathcal{A},$ $\left.
\left. \pi \right\vert _{b(\mathcal{A)}}\left( a\right) \right\vert _{%
\mathcal{D}_{\mathcal{E}}}=\pi \left( a\right) .$

\item $\pi \left( \mathcal{A}\right) ^{^{\prime }}\cap C^{\ast }(\mathcal{D}%
_{\mathcal{E}})$ \ is identified with a\textbf{\ }von Neumann algebra on $%
\mathcal{H}$. Indeed, 
\begin{eqnarray*}
\pi \left( \mathcal{A}\right) ^{^{\prime }}\cap C^{\ast }(\mathcal{D}_{%
\mathcal{E}}) &=&\{T\in B\left( \mathcal{H}\right) \cap C^{\ast }(\mathcal{D}%
_{\mathcal{E}});T\pi \left( a\right) \subseteq \pi \left( a\right) T\text{
for all }a\in \mathcal{A}\} \\
&=&\{S\in b(C^{\ast }(\mathcal{D}_{\mathcal{E}}));S\pi \left( a\right) =\pi
\left( a\right) S\ \text{ for all }a\in \mathcal{A}\} \\
&=&b(\pi \left( \mathcal{A}\right) ^{c})
\end{eqnarray*}%
where $\pi \left( \mathcal{A}\right) ^{c}=\{S\in C^{\ast }(\mathcal{D}_{%
\mathcal{E}});S\pi \left( a\right) =\pi \left( a\right) S$ for all $a\in 
\mathcal{A}\}$.

On the other hand, $\pi \left( \mathcal{A}\right) ^{c}$\ is a locally von
Neumann algebra \cite[p.4198]{D}, and then, $b(\pi \left( \mathcal{A}\right)
^{c})$ is identified with a von Neumann algebra on $\mathcal{H}$ \cite[%
Proposition 3.2]{D}.

\item By \cite[Proposition 3.2]{D}, $b(C^{\ast }(\mathcal{D}_{\mathcal{E}}))$
is identified with a von Neumann algebra on $\mathcal{H}$, which is
spatially isomorphic to the von Neumann algebra $\{T\in B\left( \mathcal{H}%
\right) ;P_{\iota }T=TP_{\iota },\forall \iota \in \Upsilon \}.$

\item Von Neumann algebras $\pi \left( \mathcal{A}\right) ^{^{\prime }}\cap
C^{\ast }(\mathcal{D}_{\mathcal{E}})$ and $\left. \pi \right\vert _{b(%
\mathcal{A)}}\left( b(\mathcal{A})\right) ^{\prime }\cap b(C^{\ast }(%
\mathcal{D}_{\mathcal{E}}))$ are isomorphic.
\end{enumerate}
\end{remark}

\begin{definition}
\cite[Definition 4.1]{BGK} Let $\varphi ,\psi \in \mathcal{CPCC}_{\text{loc}%
}(\mathcal{A},C^{\ast }(\mathcal{D}_{\mathcal{E}}))$. We say that $\psi $ is
dominated by $\varphi $,and note by $\varphi \geq \psi $, if $\varphi -\psi
\in \mathcal{CPCC}_{\text{loc}}(\mathcal{A},C^{\ast }(\mathcal{D}_{\mathcal{E%
}})).$
\end{definition}

\ \ \ Let $\varphi ,\psi \in \mathcal{CPCC}_{\text{loc}}(\mathcal{A},C^{\ast
}(\mathcal{D}_{\mathcal{E}}))$ such that $\psi $ is dominated by $\varphi $.
If $(\pi _{\varphi },V_{\varphi },$\ \ \ $\{\mathcal{H}^{\varphi },\mathcal{E%
}^{\varphi },\mathcal{D}_{\mathcal{E}^{\varphi }}\})$ is a minimal
Stinespring dilation associated to $\varphi $, then, by Radon Nikodym type
theorem \cite[Theorem 4.5]{BGK}, there is a unique element $T\in \pi
_{\varphi }\left( \mathcal{A}\right) ^{^{\prime }}\cap C^{\ast }(\mathcal{D}%
_{\mathcal{E}^{\varphi }})$ such that 
\begin{equation*}
\psi \left( a\right) =\varphi _{T}\left( a\right) =\left. V_{\varphi }^{\ast
}T\pi _{\varphi }\left( a\right) V_{\varphi }\right\vert _{\mathcal{D}_{%
\mathcal{E}}}
\end{equation*}%
for all $a\in \mathcal{A}.$

\begin{definition}
Let $\varphi \in \mathcal{CPCC}_{\text{loc}}(\mathcal{A},C^{\ast }(\mathcal{D%
}_{\mathcal{E}}))$. We say that $\varphi $ is pure if whenever $\psi \in 
\mathcal{CPCC}_{\text{loc}}(\mathcal{A},C^{\ast }(\mathcal{D}_{\mathcal{E}%
})) $ and $\varphi \geq \psi ,$ there is a positive number $\alpha $ such
that $\psi =\alpha \varphi .$
\end{definition}

Let $\varphi \in \mathcal{CPCC}_{\text{loc}}(\mathcal{A},C^{\ast }(\mathcal{D%
}_{\mathcal{E}}))$ and $\left( \pi _{\varphi },V_{\varphi },\{\mathcal{H}%
^{\varphi },\mathcal{E}^{\varphi },\mathcal{D}_{\mathcal{E}^{\varphi
}}\}\right) $ be a minimal Stinespring dilation associated to $\varphi $.
Then $\varphi $ is pure if and only if $\pi _{\varphi }\left( \mathcal{A}%
\right) ^{^{\prime }}\cap C^{\ast }(\mathcal{D}_{\mathcal{E}^{\varphi }})$ $%
=\{\alpha $id$_{\mathcal{D}_{\mathcal{E}^{\varphi }}};\alpha \in \mathbb{C}%
\}.$

\begin{proposition}
Let $\varphi \in \mathcal{CPCC}_{\text{loc}}(\mathcal{A},C^{\ast }(\mathcal{D%
}_{\mathcal{E}}))$. If $\varphi $ is pure, then it is a bounded operator
valued completely positive map.
\end{proposition}

\begin{proof}
Let $\left( \pi _{\varphi },V_{\varphi },\{\mathcal{H}^{\varphi },\mathcal{E}%
^{\varphi },\mathcal{D}_{\mathcal{E}^{\varphi }}\}\right) $ be a minimal
Stinespring dilation associated to $\varphi $. Since $\varphi $ is pure, $%
\pi _{\varphi }\left( \mathcal{A}\right) ^{^{\prime }}\cap C^{\ast }(%
\mathcal{D}_{\mathcal{E}^{\varphi }})=\{\alpha $id$_{\mathcal{D}_{\mathcal{E}%
}};\alpha \in \mathbb{C}\}$. On the other hand, for each $\iota \in \Upsilon
,$ $\left. P_{\iota }\right\vert _{\mathcal{D}_{\mathcal{E}^{\varphi }}}\in
\pi _{\varphi }\left( \mathcal{A}\right) ^{^{\prime }}\cap C^{\ast }(%
\mathcal{D}_{\mathcal{E}^{\varphi }})$, and so $P_{\iota }=$id$_{\mathcal{H}%
^{\varphi }}$. Therefore, for each $\iota \in \Upsilon ,$ $\mathcal{H}%
_{\iota }^{\varphi }=\mathcal{H}^{\varphi }$, $\mathcal{D}_{\mathcal{E}%
^{\varphi }}=\mathcal{H}^{\varphi }$ and $C^{\ast }(\mathcal{D}_{\mathcal{E}%
^{\varphi }})=B\left( \mathcal{H}^{\varphi }\right) $. Consequently, $\pi
_{\varphi }\left( a\right) \in B\left( \mathcal{H}^{\varphi }\right) $ for
all $a\in \mathcal{A}.\ $Therefore, $V_{\varphi }^{\ast }\pi _{\varphi
}\left( a\right) V_{\varphi }\in B\left( \mathcal{H}\right) $ for all $a\in 
\mathcal{A}$, and 
\begin{equation*}
\varphi \left( a\right) =\left. V_{\varphi }^{\ast }\pi _{\varphi }\left(
a\right) V_{\varphi }\right\vert _{\mathcal{D}_{\mathcal{E}}}\in b\left(
C^{\ast }(\mathcal{D}_{\mathcal{E}})\right)
\end{equation*}%
for all $a\in \mathcal{A}$.
\end{proof}

\begin{theorem}
\label{8} Let $\mathcal{A}$ and $\mathcal{B}$ be two unital locally $C^{\ast
}$-algebras, $\varphi \in \mathcal{CPCC}_{\text{loc}}(\mathcal{A},\mathcal{D}%
_{\mathcal{E}}),$ $\psi \in \mathcal{CPCC}_{\text{loc}}(\mathcal{B},C^{\ast
}(\mathcal{D}_{\mathcal{F}}))$ and $\phi \in \mathcal{CPCC}_{\text{loc}}(%
\mathcal{A\otimes B},C^{\ast }(\mathcal{D}_{\mathcal{E\otimes F}}))$. If $%
\phi $ is dominated by $\varphi \otimes \psi $ and $\varphi $ is pure, then
there exists $\widetilde{\psi }\in \mathcal{CPCC}_{\text{loc}}(\mathcal{B}%
,C^{\ast }(\mathcal{D}_{\mathcal{F}}))$ dominated by $\psi $ such that $\phi
=\varphi \otimes \widetilde{\psi }.$
\end{theorem}

\begin{proof}
Let $\left( \pi _{\varphi },V_{\varphi },\{\mathcal{H}^{\varphi },\mathcal{E}%
^{\varphi },\mathcal{D}_{\mathcal{E}^{\varphi }}\}\right) $ and $\left( \pi
_{\psi },V_{\psi },\{\mathcal{H}^{\psi },\mathcal{E}^{\psi },\mathcal{D}_{%
\mathcal{E}^{\psi }}\}\right) $ be the minimal Stinespring dilations
associated to $\varphi $ and $\psi $. Since $\varphi $ is pure, $\mathcal{D}%
_{\mathcal{E}^{\varphi }}=\mathcal{H}^{\varphi }$ and $C^{\ast }(\mathcal{D}%
_{\mathcal{E}^{\varphi }})=B\left( \mathcal{H}^{\varphi }\right) $. By
Proposition \ref{3}, $\ (\pi _{\varphi }\otimes \pi _{\psi },V_{\varphi
}\otimes V_{\psi },\{\mathcal{H}^{\varphi }\otimes \mathcal{H}^{\psi },%
\mathcal{H}^{\varphi }\otimes \mathcal{E}^{\psi },$ $\mathcal{D}_{\mathcal{H}%
^{\varphi }\otimes \mathcal{E}^{\psi }}\})$ is a minimal Stinespring
dilation associated to $\varphi \otimes \psi $. We have: 
\begin{eqnarray*}
\left( \left( \pi _{\varphi }\otimes \pi _{\psi }\right) \left( \mathcal{%
A\otimes B}\right) \right) ^{\prime } &=&\left( \left( \left. \pi _{\varphi
}\right\vert _{b\mathcal{(A)}}\otimes \left. \pi _{\psi }\right\vert _{b%
\mathcal{(B)}}\right) \left( b(\mathcal{A)\otimes }_{\text{alg}}b\mathcal{(B)%
}\right) \right) ^{\prime } \\
&=&\left( \left. \pi _{\varphi }\right\vert _{b\mathcal{(A)}}\left( b(%
\mathcal{A)}\right) \right) ^{\prime }\overline{\otimes }\left( \left. \pi
_{\psi }\right\vert _{b\mathcal{(B)}}\left( b\mathcal{(B)}\right) \right)
^{\prime }
\end{eqnarray*}%
where "$\overline{\otimes }$ " denotes the tensor product of von Neumann
algebras, and 
\begin{eqnarray*}
b(C^{\ast }(\mathcal{D}_{\mathcal{H}^{\varphi }\otimes \mathcal{E}^{\psi
}})) &=&\{R\in B(\mathcal{H}^{\varphi }\otimes \mathcal{H}^{\psi });R(\text{%
id}_{\mathcal{H}^{\varphi }}\otimes P_{\iota })=(\text{id}_{\mathcal{H}%
^{\varphi }}\otimes P_{\iota })R,\forall \iota \in \Upsilon \} \\
&=&B(\mathcal{H}^{\varphi })\overline{\otimes }\{S\in B(\mathcal{H}^{\psi
});SP_{\iota }=P_{\iota }S,\forall \iota \in \Upsilon \} \\
&=&B(\mathcal{H}^{\varphi })\overline{\otimes }b(C^{\ast }(\mathcal{D}_{%
\mathcal{E}^{\psi }})).
\end{eqnarray*}%
Thus, 
\begin{eqnarray*}
&&\left( \left( \pi _{\varphi }\otimes \pi _{\psi }\right) \left( \mathcal{%
A\otimes B}\right) \right) ^{\prime }\cap b(C^{\ast }(\mathcal{D}_{\mathcal{H%
}^{\varphi }\otimes \mathcal{E}^{\psi }})) \\
&=&\left( \left( \left. \pi _{\varphi }\right\vert _{b\mathcal{(A)}}\left( b(%
\mathcal{A)}\right) \right) ^{\prime }\overline{\otimes }\left( \left. \pi
_{\psi }\right\vert _{b\mathcal{(B)}}\left( b\mathcal{(B)}\right) \right)
^{\prime }\right) \cap \left( B(\mathcal{H}^{\varphi })\overline{\otimes }%
b(C^{\ast }(\mathcal{D}_{\mathcal{E}^{\psi }}))\right) \\
&=&\left( \pi _{\varphi }\left( \mathcal{A}\right) ^{\prime }\cap B(\mathcal{%
H}^{\varphi })\right) \overline{\otimes }\left( \left( \left. \pi _{\psi
}\right\vert _{b\mathcal{(B)}}\left( b\mathcal{(B)}\right) \right) ^{\prime
}\cap b(C^{\ast }(\mathcal{D}_{\mathcal{E}^{\psi }}))\right) \\
&&\text{(since }\varphi \text{ is pure) } \\
&=&\{\alpha \text{id}_{\mathcal{H}^{\varphi }};\alpha \in \mathbb{C}%
\}\otimes \left( \pi _{\psi }\left( \mathcal{B}\right) ^{\prime }\cap
C^{\ast }(\mathcal{D}_{\mathcal{E}^{\psi }})\right) .
\end{eqnarray*}%
Therefore, 
\begin{equation*}
\left( \left( \pi _{\varphi }\otimes \pi _{\psi }\right) \left( \mathcal{%
A\otimes B}\right) \right) ^{\prime }\cap C^{\ast }(\mathcal{D}_{\mathcal{H}%
^{\varphi }\otimes \mathcal{E}^{\psi }})=\{\alpha \text{id}_{\mathcal{H}%
^{\varphi }};\alpha \in \mathbb{C}\}\otimes \left( \pi _{\psi }\left( 
\mathcal{B}\right) ^{\prime }\cap C^{\ast }(\mathcal{D}_{\mathcal{E}^{\psi
}})\right) .
\end{equation*}%
Since $\phi \ $and $\varphi \otimes \psi $ are local completely contractive
and local completely positive and $\phi $ is dominated by $\varphi \otimes
\psi $, by Radon Nikodym theorem \cite[Theorem 4.5]{BGK}, there is a unique
positive contractive linear operator $R\in \left( \left( \pi _{\varphi
}\otimes \pi _{\psi }\right) \left( \mathcal{A\otimes B}\right) \right)
^{\prime }\cap C^{\ast }(\mathcal{D}_{\mathcal{H}^{\varphi }\otimes \mathcal{%
E}^{\psi }})$ such that $\phi =\left( \varphi \otimes \psi \right) _{R}$.
Therefore, there is $T\in \pi _{\psi }\left( \mathcal{B}\right) ^{\prime
}\cap C^{\ast }(\mathcal{D}_{\mathcal{E}^{\psi }})$ such that $R=$id$_{%
\mathcal{H}^{\varphi }}\otimes T$ and 
\begin{eqnarray*}
\phi \left( a\otimes b\right) &=&\left( \varphi \otimes \psi \right) _{\text{%
id}_{\mathcal{H}^{\varphi }}\otimes T}\left( a\otimes b\right) \\
&=&\left( V_{\varphi }^{\ast }\otimes V_{\psi }^{\ast }\right) \left( \text{%
id}_{\mathcal{H}^{\varphi }}\otimes T\right) \left( \pi _{\varphi }\otimes
\pi _{\psi }\right) \left( a\otimes b\right) \left( V_{\varphi }\otimes
V_{\psi }\right) \\
&=&V_{\varphi }^{\ast }\pi _{\varphi }\left( a\right) V_{\varphi }\otimes
V_{\psi }^{\ast }T\pi _{\psi }\left( b\right) V_{\psi }=\varphi \left(
a\right) \otimes \psi _{T}\left( b\right)
\end{eqnarray*}%
for all $a\in \mathcal{A}$ and $b\in \mathcal{B}$.$\ $Hence, there is $%
\widetilde{\psi }=\psi _{T}\in \mathcal{CPCC}_{\text{loc}}(\mathcal{B}%
,C^{\ast }(\mathcal{D}_{\mathcal{F}}))$ such that $\phi =\varphi \otimes 
\widetilde{\psi }$.\ Moreover, $\widetilde{\psi }$ is dominated by $\psi $.
\end{proof}

\begin{corollary}
Let $\mathcal{A}$ be a unital $C^{\ast }$-algebra, $\mathcal{B}$ be a unital
locally $C^{\ast }$-algebra, $\varphi \in \mathcal{CP}(\mathcal{A},B(%
\mathcal{H})),$ $\psi \in \mathcal{CPCC}_{\text{loc}}(\mathcal{B},C^{\ast }(%
\mathcal{D}_{\mathcal{F}}))$ and $\phi \in \mathcal{CPCC}_{\text{loc}}(%
\mathcal{A\otimes B},C^{\ast }(\mathcal{D}_{\mathcal{H\otimes F}}))$. If $%
\phi $ is dominated by $\varphi \otimes \psi $ and $\varphi $ is pure, then
there is $\widetilde{\psi }\in \mathcal{CPCC}_{\text{loc}}(\mathcal{B}%
,C^{\ast }(\mathcal{D}_{\mathcal{F}}))$ which is dominated by $\psi $ and
such that $\phi =\varphi \otimes \widetilde{\psi }$.
\end{corollary}

\begin{acknowledgement}
I wish to thank A. Dosiev and P. Santhosh Kumar for useful discussions on
the details of the present work. Also, I would like to thank the referee for
his/her careful reading and useful comments.
\end{acknowledgement}


\begin{thebibliography}{99}
\bibitem{BGK} B. V. R. Bhat, A. Ghatak and P. S. Kumar, \textit{%
Stinespring's theorem for unbounded operator valued local completely
positive maps and its applications, }Indag. Math. 32(2021),2, 547-578.

\bibitem{BO} B. V. R. Bhat and H. Osaka, \textit{A factorization property of
positive maps on }$C^{\ast }$\textit{-algebras, }Int. J. Quantum Inf.
18(2020), 5.

\bibitem{D1} A. Dosiev, \textit{Local operator spaces, unbounded operators
and multinormed }$C^{\ast }$\textit{-algebras, }J. Funct. Anal. 255(2008),
1724--1760.

\bibitem{D} A. Dosiev, \textit{Multinormed }$\mathit{W}^{\mathit{\ast }}$%
\textit{-algebra and unbounded operators}, Proc. AMS, 140 (2012) 4187-4202.

\bibitem{F} M. Fragoulopoulou, Topological algebras with involution,
Elsevier, 2005.

\bibitem{G} A. Gheondea, \textit{Operator models for Hilbert locally }$%
C^{\ast }$\textit{-modules,} Oper. Matrices 11 (2017), no. 3, 639--667.

\bibitem{I} A. Inoue, \textit{Locally }$C^{\ast }$\textit{-algebras}, Mem.
Fac. Sci. Kyushu Univ. Ser. A, \textbf{25} (1971), 197--235.

\bibitem{J} M. Joi\c{t}a, \textit{Strict completely positive linear maps
between locally }$C^{\ast }$\textit{-algebras and representations on Hilbert
modules}, J. London Math. Soc. (2), 66 (2002), no.2, pp.421-432.

\bibitem{MJ1} M. Joi\c{t}a, \textit{Unbounded local completely positive maps
of local order zero}, Positivity, (2021).

\bibitem{MJ2} M. Joi\c{t}a, \textit{Completely positive maps of order zero
on pro-}$C^{\ast }$\textit{-algebras, }Forum Math. 33(2021),1, 29-37.

\bibitem{MJ3} M. Joi\c{t}a, \textit{Decomposable local positive maps on
locally }$C^{\ast }$\textit{-algebras, Manuscript.}

\bibitem{Ph} N.C.\ Phillips, \textit{Inverse limits of }$C^{\ast }$\textit{%
-algebras}, J. Operator Theory 19 (1988), no. 1, 159--195.
\end{thebibliography}
\end{document}